\documentclass[11pt,reqno]{amsart}
\usepackage{amsmath, euscript, amssymb, amsfonts}

\theoremstyle{plain}
\newtheorem{theorem}{Theorem}
\newtheorem{lemma}{Lemma}
\newtheorem{corollary}{Corollary}

\theoremstyle{definition}
\newtheorem{definition}{Definition}

\newtheorem*{acknowledgements}{Acknowledgements}

\theoremstyle{remark}
\newtheorem{remark}{Remark}

\newcommand{\oR}{{\mathbb R}}

\newcommand{\oC}{{\mathbb C}}
\newcommand{\oZ}{{\mathbb Z}}

\def\d{{\mathrm d}}
\newcommand{\pr}{\mathop{\rm pr}\nolimits}
\newcommand{\supp}{\mathop{\rm supp}\nolimits}

\newcommand{\R}{\mathop{\rm Re}\nolimits}
\newcommand{\I}{\mathop{\rm Im}\nolimits}
\newcommand{\Int}{\mathop{\rm Int}\nolimits}
\newcommand{\ch}{\mathop{\rm ch}\nolimits}

\begin{document}

 \title[Decomposition theorems and kernel theorems]{Decomposition theorems
  and kernel theorems \\for a class
of functional spaces}

\author{M.~A.~Soloviev}
\address{Lebedev
Physical Institute, Russian Academy of Sciences,
 Leninsky
Prospect 53, Moscow 119991, Russia} \email{soloviev@lpi.ru}
\thanks{}

 \begin{abstract}
We prove new theorems about properties of  generalized functions
defined on  Gelfand-Shilov spaces $S^\beta$ with $0\le\beta<1$.
For each open cone $U\subset \oR^d$ we define a space $S^\beta(U)$
which is related to $S^\beta(\oR^d)$ and consists of  entire
analytic functions rapidly decreasing inside  $U$ and having order
of growth $\le 1/(1-\beta)$ outside the cone. Such sheaves of
spaces arise naturally in nonlocal quantum field theory, and this
motivates our investigation. We prove that the spaces $S^\beta(U)$
are complete and nuclear and establish a decomposition theorem
which implies that every continuous functional defined on
$S^\beta(\oR^d)$ has a unique minimal closed carrier cone in
$\oR^d$. We also prove kernel theorems for  spaces over open and
closed cones and elucidate the relation between the carrier cones
of multilinear forms and those of the generalized functions
determined by these forms.
  \end{abstract}
\maketitle

 \section{Introduction}
In this paper we investigate the  angular localizability property
of  generalized functions belonging to the spaces $S^{\prime
\beta}$,  $0\le\beta<1$. This property was revealed in
applications of these  classes of generalized  functions to
nonlocal quantum field theory. The test function spaces $S^\beta$
and $S^\beta_\alpha$ were introduced by Gelfand and
Shilov~\cite{GS}. If $\beta<1$, they consist of entire functions,
and   continuous linear functionals on these spaces are analytic
in the sense that they  are representable by Taylor series
convergent in the  topology (weak as well as strong) of the dual
spaces $S^{\prime\,\beta}$ and $S^{\prime\,\beta}_\alpha$.

 To develop a nonlocal field theory, one should  use certain spaces
 related to $S^\beta$, $S^\beta_\alpha$ and associated with
 cones in $\oR^d$. The notion of a minimal carrier
cone of an analytical functional plays a key role in these
applications.  The existence of such a quasi-support follows from
appropriate  decomposition theorems for  spaces over cones. For
spaces with two indices,  the corresponding analysis is performed
in~\cite{S1}. It is based on the fact that  these spaces belong to
the well-studied class of DFS spaces, that is, spaces  dual to FS
(Fr\'echet-Schwartz) spaces. The properties of DFS spaces are
reviewed, for example, in the  survey~\cite{Zh}. The topological
structure of the spaces $S^\beta(U)$, where $U$ is an open cone in
$\oR^d$, is more complicated than that of $S^\beta_\alpha(U)$.
They are not DFS spaces, and even the proof of their completeness
is a challenge. If $\beta=0$, then proving decomposition theorems
presents additional difficulties, and one way around them is to
use methods of the theory of hyperfunctions. This limiting case is
of prime interest because the space  $S^0$ is nothing but the
Fourier transform of the Schwartz space $\mathcal D$ of all
infinitely differentiable  compactly supported functions. The
existence of smallest carrier cones for elements of the dual space
$S^{\prime 0}$ can be established in a roundabout way~\cite{S} by
restricting  these functionals to $S^0_\alpha$. But this result
alone is not sufficient for  applications, and the properties of
these function spaces call for further  investigation. In this
paper, particular attention is given to the   extension of the
theory  to multilinear forms, including the derivation of  kernel
theorems for  $S^\beta(U)$ and for spaces over closed cones, which
are constructed from  spaces over open cones by means of the
inductive limit.

We note that spaces of type $S^\beta$ over cones arise naturally
when one extends the theory of  Fourier-Laplace transformation to
analytic functionals, and especially when one  generalizes the
Paley-Wiener-Schwartz theorem~\cite{V,H1}. Vladimirov's
version~\cite{V} of this theorem shows that there is an
isomorphism between the space of tempered distributions supported
in a  properly convex closed cone $K$ and an algebra of analytic
functions defined on a certain tubular domain and growing at most
polynomially. If we relax the bound on the growth of functions at
infinity and pose the question of finding the class of functionals
corresponding to the enlarged algebra, then we inevitably arrive
at the spaces $S^{\prime\beta}(K)$. If  the bound is removed
altogether, then we arrive at $S^{\prime 0}(K)$.

The paper is organized as follows. The next section contains basic
definitions and  some  preliminary information on the function
spaces in question. In Sect.~3 we use Palomodov's
criterion~\cite{P} to prove that these spaces are complete. In
Sects.~4, 5 we prove the main decomposition theorem for spaces
over cones. This theorem implies that the correspondence $K\to
S^\beta(K)$ is a lattice (anti-)homomorphism. For $\beta>0$, the
proof is simpler and uses the non-triviality of the space
$S^\beta_{1-\beta}$. The proof for  $\beta=0$ is given in Sect.~5
and relies on  H\"ormander's estimates~\cite{H2}. Since the weight
functions in these estimates are plurisubharmonic, we develop a
general technique (Appendix 1) for approximating the indicator
functions (that determine the spaces under study) by
plurisubharmonic functions. In Sect.~6 we prove that every element
of $S^{\prime \beta}(\oR^d)$ has a unique minimal closed carrier
cone. In Sect.~7 we establish that the spaces associated with open
and closed cones are nuclear and indicate some consequences of
this result. The corresponding kernel theorems, which enable one
to identify   multilinear separately continuous forms on these
spaces with linear functionals, are proved in Sect.~8. The method
devised for this purpose is also applicable to other spaces of
analytic functions. In Sect.~9 we show that the study of analytic
functionals generated by multilinear forms leads naturally to the
notion of a strong carrier cone. The difference between the
notions of a carrier cone and a strong carrier cone is elucidated
in Appendix~2. In Sect.~10 we derive a Paley-Wiener-Schwartz-type
theorem, which  precisely describes the properties of the Laplace
transforms of  functionals that belong to $S^{\prime \beta}$ and
are strongly carried by a properly convex cone.

The theorems of  the theory of linear topological spaces used
below are contained in~\cite{G,RR,Sch}. We refer to~\cite{V} for
the basic facts about plurisubharmonic functions.

\section{Basic definitions and preliminaries}

Let $0\le\beta<1$ and let $U$ be an open cone in $\oR^d$. We
define  $S^{\beta,b}(U)$ to be the intersection (projective limit)
of the Hilbert spaces $H^{\beta,B}_N(U)$, $B>b$, $N=0,1,2,\dots$,
that consist of entire functions on $\oC^d$  and are equipped with
the inner products
  \begin{equation} {\langle f,\,g\rangle}_{U,B,N} = \int
\overline{f(z)}
g(z)\prod_{j=1}^d(1+|x_j|)^{2N}\exp\{-2d(Bx,U)^{1/(1-\beta)} -
2|By|^{1/(1-\beta)}\}\,d\lambda ,
   \label{1*}
   \end{equation}
   where $z=x+iy$,
$d(x,U)=\inf_{\xi\in U}|x-\xi|$ is the distance from $x$ to $U$
and $d\lambda=dx\,dy$ is the Lebesgue  measure on $\oC^n$. Note
that $d(Bx,U)=Bd(x,U)$ because $U$ is a cone. It is clear from
this definition that $S^{\beta,b}(U)$ is a Fr\'echet space (that
is, a complete  metrizable space). We denote  the union of
$S^{\beta,b}(U)$ over all $b>0$ by $S^\beta(U)$ and endow it with
the inductive limit topology.  This space is independent of the
choice of the norm $|\cdot|$ on $\oR^d$, because all these norms
are equivalent. In what follows we use either the Euclidean norm
or the $l^{1/(1-\beta)}$--norm
 \begin{equation}
\left(\sum_{j=1}^d|y_j|^{1/(1-\beta)}\right)^{1-\beta}.
\label{norm}
   \end{equation}
The latter is convenient when treating miltilinear forms on
$S^0(\oR^{d_1})\times\dots\times S^0(\oR^{d_n})$ because then the
weight function  in~\eqref{1*} has a  multiplicative property.
   Namely, if we put
   \begin{equation}
\rho_{U,B,N}(z)=-N\sum_j\ln(1+|x_j|)+d(Bx,U)^{1/(1-\beta)}
+|By|^{1/(1-\beta)},
    \label{rho}
   \end{equation}
then the function determining  the space $S^0(U_1\times\dots\times
U_n)$, where $U_k\subset \oR^{d_k}$, is given by the product of
the functions of the functions determining the $S^0(U_k)$:
\begin{equation}
\exp\{-\rho_{U_1\times\dots\times U_n,B,N}(z)\} =\prod_{k=1}^n
\exp\{-\rho_{U_k,B,N}(z_k)\}.
 \label{m}
  \end{equation}

The space  $S^{\beta,b}(U)$ may also be  represented as the
intersection of the Banach spaces    $E^{\beta,B}_N(U)$ of entire
functions with the norms
    \begin{equation}
     \|f\|'_{U,B,N}= \sup_{z\in \oC^d}|f(z)|\,e^{-\rho_{U,B,N}(z)}.
    \label{2*}
    \end{equation}
This is precisely the  original definition given in~\cite{FS}, but
the reformulation in terms of Hilbert spaces  is best suited to
most of  the questions discussed below. The representation
$$
  S^{\beta, b}(U)=\bigcap_{B>b, N\ge 0}E^{\beta,B}_N(U)
$$
makes it clear that $S^\beta(\oR^d)$ coincides with the
Gelfand-Shilov space  $S^\beta$ and with the Gurevich space
$W^\Omega$, where $\Omega(y)=y^{1/(1-\beta)}$. However, these
spaces were not given a topology in~\cite{GS} and the  notion of
convergence of sequences was used instead. In Sect.~3 we show that
this simplified ``sequential'' approach  agrees with the natural
topology described above. We also note that $S^{\beta, b}(U)$ can
be treated as a countably normed space $Z(M_p)$ specified by
$M_p=\exp\{-\rho_{U,b+1/p,p}\}$ if we omit the condition
$M_p(z)\ge C(y)$  from the definition~\cite{GS} of this class of
spaces. The equivalence between  the system~\eqref{2*} of norms
and the system $\|f\|_{U,B,N}$ defined by the inner
products~\eqref{1*} can be established using Cauchy's integral
formula, which shows that
\begin{equation}
|f(z)|\le C\|f\|_{L^2(\mathcal B)},
 \notag
  \end{equation}
where $\mathcal B$ is any bounded neighborhood of  $z$ in $\oC^d$.
Taking  $\mathcal B=\{\zeta\colon |z-\zeta|<1\}$ and applying the
triangle   inequality to every term on the right-hand side
of~\eqref{rho},   we see that $\rho_{U,B',N}(\zeta)\le
\rho_{U,B,N}(z)+C_{B,B',N}$ for $\zeta\in \mathcal B$ and  any
$B'<B$. Therefore,
\begin{equation}
|f(z)|^2 e^{-2\rho_{{}_{U, B,N}}(z)}\le C'\int_{\mathcal
B}|f(\zeta)|^2 e^{-2\rho_{{}_{U, B',N}}(\zeta)}\d \lambda \le C'
\|f\|^2_{U, B', N},
 \label{3*}
  \end{equation}
On the other hand, it is clear that $\|f\|_{U,B,N}\le C''
\|f\|'_{U,B',N+d+1}$.

For each closed cone $K\subset\oR^d$ we define the space
$S^\beta(K)$  as the inductive limit of the spaces $S^\beta(U)$,
where $U$ runs through the open cones that contain $K$ as a
compact subcone. (This  is written\footnote{For  arbitrary cones
$V_1$, $V_2$, the notation  $V_1\Subset V_2$ means that $\Bar
V_1\setminus \{0\}\subset V_2$. Here and in what follows, we use a
bar to denote the closure of a set.}  $K\Subset U$.) All these
spaces are continuously embedded into the space $S^\beta(\{0\})$
associated with the degenerate closed cone consisting of one
point, namely,  the origin. Its elements are entire functions of
order $1/(1-\beta)$  and finite type or of order less than
$1/(1-\beta)$. It should be noted that we suggest cones for
geometric visualization, although we are really   dealing with a
sheaf of spaces over the sphere  compactifying $\oR^d$. Although
the cone $\{0\}$ is closed in $\oR^d$, it corresponds to the empty
subset of the sphere, which is both closed and open. Therefore the
space $S^\beta(\{0\})$ along with  its topology can be defined
directly by formula~\eqref{1*} with $d(x,0)=|x|$. As we shall see
in Sect.~5, this topology coincides with the inductive topology
determined by the injections $S^\beta(U)\to S^\beta(\{0\})$, where
$U$ ranges over all open cones in $\oR^d$.

\begin{definition} Let $v$ be a continuous linear  functional on
$S^\beta(\oR^d)$. We say that $v$ is {\it carried} by a closed
cone  $K\subset \oR^d$  if this functional admits a continuous
extension to $S^\beta(K)$.
\end{definition}
When such an extension exists, it is unique by  the following
theorem.

 \begin{theorem} There is a constant $\lambda$
$($depending only on $d$ and $\beta$$)$ such that the space
$S^{\beta,\lambda b}(\oR^d)$ is dense in $S^{\beta,b}(U)$ in the
topology of $S^{\beta,\lambda b}(U)$  for every open cone
$U\subset\oR^d$ and  every $b>0$. As a consequence,
$S^\beta(\oR^d)$ is sequentially dense in $S^\beta(U)$ and in any
space $S^\beta(K)$, where $K$ is a closed cone.
 \end{theorem}
\begin{proof} When $\beta>0$, we can use the fact that the space
$S^\beta_{1-\beta}$ is nontrivial. According to~\cite{GS},  there
is a $\gamma>0$ such that for any $A>0$ the space
$S^\beta_{1-\beta}(\oR^d)$ contains a nontrivial nonnegative
function $g_0$ satisfying the bound
  \begin{equation}
|g_0(z)|\le C \exp\left\{-\left|\frac{x}{A}\right|^{1/(1-\beta)} +
\left|\frac{\gamma y}{A}\right|^{1/(1-\beta)}\right\}.
 \label{g}
  \end{equation}
  Let $f\in S^{\beta,b}(U)$. We normalize $g_0$ by the condition
$\int g_0(x)\,d x=1$ and set $f_\nu=f \sigma_\nu$, where
$\sigma_\nu(z)$ is a sequence of Riemann sums for the integral
$\int g_0(z-\xi)\, d\xi$ or, more explicitly,
$$
\sigma_\nu (z) = \sum_{k\in \oZ^n,|k|<\nu^2}g_0
\left(z-\frac{k}{\nu}\right)\nu^{-n}.
$$
 Clearly, $f_\nu \in S^\beta(\oR^d)$ if $A< 1/b$.
 The sequence $\sigma_\nu$ converges to $1$ in $\oR^d$ and we have
 $|\sigma_\nu (z)| \le C'\exp\{\,|\gamma y/A|^{1/(1-\beta)}\}$
 because the integral sums for  $\exp\{-|x/A|^{1/(1-\beta)}\}$
 are bounded.   Therefore, $\sigma_\nu(z)\to 1$ uniformly on
 compact sets in $\oC^d$  by  Vitali's theorem. Moreover,
 the sequence $f_\nu$ is bounded in every norm
$|\cdot|'_{U,B,N}$, where $B>b+\gamma/A$. Hence $f_\nu\to f$ in
the topology of $S^{b+\gamma/A}(U)$ because it is a Montel
space.\footnote{In~\cite{GS},  countably normed Montel spaces were
termed {\it perfect spaces} and it was  proved  that every space
$Z(M_p)$ is perfect. We also note that $S^b(U)$ is a Montel space
because it is nuclear, see below.} We see that in this case
Theorem 1 holds for any $\lambda>1+\gamma$.

If $\beta=0$, then this argument fails because $S^0_1$ is trivial.
However, Theorem 1 can easily be deduced from an analogous theorem
established for $S^0_\alpha(U)$ in~\cite{S1} by  an alternative
method using H\"ormander's $L^2$-estimates. The space
$S^0_\alpha(U)$ with $\alpha>1$ is the union of the Banach spaces
$E^{0,B}_{\alpha,A}(U)$ of entire functions with the norms
\begin{equation}
\|f\|_{\alpha, U,A,B}= \sup_{z\in
\oC^d}|f(z)|\,\exp\left\{\left|\frac{x}{A}\right|^{1/\alpha}-
d(Bx,U) -|By|\right\}.
 \label{4*}
   \end{equation}
 Let us show that $S^0_\alpha(U)$ is dense in $S^0(U)$.  Let $f\in
S^{0,b}(U)$. The sets
$$\{\tilde f\colon \|f - \tilde f\|'_{U,B,N}
\le \delta\},\quad \text{where}\,\, \delta>0,\, B>b,\,
N=0,1,\dots,
$$
 form a base of neighborhoods of $f$ in $S^{0,b}(U)$.
We take a function $g \in E^{0,1}_{\alpha,1}(\oR^d)$ with
$g(0)=1$ and consider the sequence $f_\nu(z)=f(z)g(z/\nu)$. If
$\nu>1/(B-b)$, then $f_\nu\in E^{0,B}_{\alpha,\nu}(U)$.
 Let $b<B'<B$. Then we have
$$
\|f - f_\nu\|'_{U,B,N} \le \|f\|'_{U,B',N+1} \sup_z
\left|1-g\left(\frac{z}{\nu}\right)\right|(1+|x|)^{-1}
e^{-\epsilon|y|},
$$
where $\epsilon=B-B'>0$. If $\nu>2/\epsilon$, then
$|1-g(z/\nu)|\le C \exp\{\epsilon|y|/2\}$. Given $\delta>0$, we
can choose  $R$ such that
$$
C\|f\|'_{U,B',N+1}\sup_{|z|\ge R}(1+|x|)^{-1} e^{-\epsilon|y|/2}<
\delta.
$$
 Taking    $\nu$ large enough for the inequality
$\|f\|'_{U,B',N+1}\sup_{|z|<R} |1-g(z/\nu)|< \delta$ to hold, we
obtain $\|f-f_\nu\|'_{U,B,N}<\delta$. Now let $\alpha'>\alpha$.
Clearly, $f_\nu\in E^{0,B}_{\alpha',A}(U)$ for any $A>0$.
Theorem~5 of~\cite{S1} shows that if $\mu>2ed$, then the function
 $f_\nu$ can be approximated in the norm $\|\cdot\|_{\alpha',
U,1,\mu B}$ by elements of $E^{0,\mu B}_{\alpha',1}(\oR^d)$ with
any degree of accuracy. This norm is stronger than
$\|\cdot\|_{U,\mu B,N}$. Therefore, in this case Theorem 1 holds
for any $\lambda>2ed$, because then there are $\mu>2ed$ and $B>b$
such that $\mu B<\lambda b$.
\end{proof}
 We note that if $v\in
S^{\prime\, \beta}$ is carried by a cone $K$, then the restriction
of $v$ to $S^{\beta_1}$ with $\beta_1<\beta$ is also carried by
this cone. However the converse is not  true in general.  In what
follows, the open mapping theorem is  used repeatedly. When
dealing with the $S^\beta$-type spaces associated with cones, we
can use Grothendieck's version~\cite{G} of this important theorem
(or the even more general version given  by Raikov in Appendix~1
to the Russian edition of~\cite{RR}) because all these spaces,
being Hausdorff and inductive limits of sequences of Fr\'echet
spaces, belong to the class\footnote{We note that, in contrast to
the definition given in ~\cite{Sch}, Grothendieck's definition of
this class does not require that the inductive limit be strict.}
$\mathcal{LF}$ and  are  ultrabornological (spaces of type
$(\beta)$ in the terminology of~\cite{G}). In~\cite{S5}, we showed
that none of the spaces $S^\beta(U)$, $S^\beta(K)$ except
$S^\beta(\{0\})$ are  DFS spaces because their dual spaces are
non-metrizable.

 \section{The completeness theorem}
The completeness of  the spaces $S^\beta(U)$ was proved
in~\cite{S3} using another definition  given in terms of real
variables. For the reader's convenience, we present an alternative
proof starting from the norms~\eqref{2*}.
\begin{theorem} The inductive spectrum of the Fr\'echet
spaces $S^{\beta,b}(U)$, $b=1,2,\dots$, is acyclic.
\end{theorem}
\begin{proof} Let  ${\mathcal U}_{b}$ be the neighborhood of the
origin in $S^{\beta,b}(U)$ specified by $\|f\|'_{U,b+1/2,0}<1/2$.
Clearly,  ${\mathcal U}_{b_0}\subset {\mathcal U}_{b}$ for any
$b>b_0$. According to Theorem 6.1 of~\cite{P}, it suffices to
verify that the topology on   ${\mathcal U}_{b_0}$ induced by that
of $S^{\beta,b}(U)$,   $b>b_0$,   is independent of $b$. Let
$f_0\in {\mathcal U}_{b_0}$ and $B>b$.   We denote by ${\mathcal
V}_{B,N,\epsilon}$  the intersection of ${\mathcal U}_{b_0}$ and
the neighborhood of  $f_0$ in $S^{\beta,b}(U)$ given by
$\|f-f_0\|_{U,B,N}<\epsilon$. We shall show that for any $B$,
$B_1$ satisfying $B>B_1>B_0=b_0+1/2$, for every $N_1\ge0$, and for
every $\epsilon_1>0$, there are numbers $N$ and $\epsilon$ such
that
\begin{equation}
 {\mathcal V}_{B,N,\epsilon}\subset {\mathcal V}_{B_1,N_1,\epsilon_1}.
 \label{x0}
   \end{equation}
This means that the topology induced on ${\mathcal U}_{b_0}$ by
that of $S^{\beta,b}(U)$ is not weaker than the topology induced
by that of  $S^{\beta,b_1}(U)$, where $b_1<b$ (the reverse is
obvious). In what follows we set $\beta=0$ for simplicity and
comment on the case $\beta>0$ at the end of the proof. If $f\in
{\mathcal V}_{B,N,\epsilon}$,  then the function  $f_1=f-f_0$
satisfies the  estimates
     \begin{gather}
     |f_1(x)|<e^{B_0d(x,U) +B_0|y|},
    \label{x1}\\
 |f_1(x)|<\epsilon\,
   (1+|x|)^{-N}e^{Bd(x,U) +B|y|}.
   \label{x2}
    \end{gather}
   We claim that, for properly chosen $N$ and $\epsilon$, this
    implies that
      \begin{equation}
    |f_1(x)|<\epsilon_1   (1+|x|)^{-N_1}e^{B_1d(x,U) +B_1|y|}.
   \label{x3}
   \end{equation}
  We introduce the notation $\varepsilon(x)=\epsilon (1+|x|)^{-N}$,
  $\varepsilon_1(x)=\epsilon_1 (1+|x|)^{-N_1}$ and define
  a number $R(x)$    by the equation
    \begin{equation}
    e^{B_0R}=\varepsilon_1 e^{B_1R}.
    \label{x4}
   \end{equation}
   In the region $d(x,U)+|y|\ge R(x)$, the inequality~\eqref{x3}
   follows from~\eqref{x1}. In the complementary region, \eqref{x3}
   follows from~\eqref{x2} if
    \begin{equation}
     \varepsilon\,e^{BR}=\varepsilon_1\,
    e^{B_1R}
    \label{x5}
   \end{equation}
   and, a fortiori,   if $\varepsilon\,e^{BR}<\varepsilon_1\, e^{B_1R}$.
   Equations~\eqref{x4} and \eqref{x5}  give $\varepsilon=\varepsilon_1^A$,
   where $A=(B-B_0)/(B_1-B_0)$. Hence the desired inclusion~\eqref{x0}
   follows if we take $\epsilon\leq \epsilon_1^A$  and $N\geq AN_1$,
   This proof extends to $\beta>0$ by an obvious
   change of notation, which yields the same conclusion with
   the modified number $A=(\tilde B- \tilde B_0)/(\tilde B_1-\tilde B_0)$,
   where     $\tilde B=B^{1/(\beta-1)}$.
\end{proof}
   By~\cite{P}, the acyclicity  ensures that the  following assertions hold.

     \begin{corollary} The space $S^{\beta}(U)$ is
     Hausdorff and complete.  A set ${\mathcal B}\subset S^{\beta}(U)$
     is bounded if and only  if it is contained in some space
     $S^{\beta,B}(U)$ and is bounded in each of its norm.
\end{corollary}
It is certainly obvious that  $S^{\beta}(U)$ is a Hausdorff space
because its topology is stronger than the topology of uniform
convergence. We also note that a linear map of $S^\beta(U)$ (as
 of any   bornological space) to a
locally convex space is continuous if and only if it is bounded on
bounded sets, which is in turn  equivalent to the sequential
continuity, see~\cite{Sch}.

 \section{Proof of the decomposition theorem for  $\beta>0$}

Here and in Sect.~5 we use the Euclidean norm on $\oR^d$. We
recall that the intersection of a cone $V$ with the unit sphere is
called the {\it projection} of this cone and is denoted by $\pr
V$.

 \begin{theorem}  Let $U$, $U_1$,  $U_2$ be open cones in
 $\oR^d$ such that $\bar U_1\cap \bar U_2\Subset U$. Then every
 function $f\in S^\beta(U)$, $\beta\ge 0$, can be decomposed as
 $f=f_1+ f_2$, where $f_i\in S^\beta(U\cup U_i)$, $i=1,2$.
\end{theorem}
\begin{proof} If $\beta>0$, then we can use the same function $g_0\in
S^\beta_{1-\beta}$ as in the proof of Theorem~1. For simplicity,
we assume that $f$ satisfies the estimate
\begin{equation}
|f(z)|\le C_N(1+|x|)^{-N}\,\exp\{d(x,U)^{1/(1-\beta)} +
|y|^{1/(1-\beta)}\},\quad N=0,1,2,\dots.
   \label{4.1}
   \end{equation}
This does not cause any loss of generality because the spaces
involved are invariant under the dilation $f(x)\to f(\lambda x)$,
$\lambda>0$. The hypothesis  $\Bar U_1 \cap\Bar U_2\Subset U$
implies that the closed cones $V_1=\Bar U_1\setminus U$, $V_2=\Bar
U_2\setminus U$ have disjoint projections. Therefore the distances
from $\pr V_1$ to $V_2$ and from $\pr V_2$ to $V_1$ are positive.
In the Euclidean metric, these distances  coincide. Indeed, if the
first distance is attained at  points $x_1\in \pr V_1$ and $x_2\in
V_2$, then the equation
$|x_1-x_2|^2=\left|\,|x_2|x_1-|x_2|^{-1}x_2\right|^2$ implies that
$d(\pr V_1, V_2)\ge d(\pr V_2, V_1)$, and the reverse inequality
holds by symmetry. We  denote this distance by $\theta$.

We now introduce the auxiliary open cone
$$
 W=\left\{\xi\in \oR^d\colon d(\xi, V_2)< \frac{\theta}{2}|\xi|\right\}
 $$
and define $g(z)$ by
\begin{equation}
 g(z)=\int\limits_W\!g_0(z-\xi)\,{\rm d}\xi .
 \label{4.2}
  \end{equation}
  We claim that if the constant $A$ in~\eqref{g} is small enough, then
  $gf\in  S^\beta(U\cup U_1)$, that is
\begin{equation}
|(gf)(z)|\le C'_N(1+|x|)^{-N}\,\exp\{d(Bx,U\cup U_1)^{1/(1-\beta)}
+ |By|^{1/(1-\beta)}\}
   \label{4.3}
   \end{equation}
 for some $B>0$.      Let $W_1=\{x\in \oR^d\colon
  d(x, V_2)\ge 3\theta|x|/4\}$. Then
\begin{equation}
|x-\xi|\ge \frac{\theta}{4}|x|\quad \text{for all}\quad x\in
W_1,\, \xi\in W.
   \label{4.4}
   \end{equation}
Indeed, if this is not the case, then there are points $x\in W_1$
and $\xi\in W$ such that $|x|=1$,  $|x-\xi|< \theta/4$, and
$|\xi|\le 1$. Also, there is  a point $x_2\in V_2$ such that
$|\xi- x_2|<\theta/2$. Then $|x-x_2|<3\theta/4$ by  the triangle
inequality. This contradicts the definition of $W_1$. It follows
from~\eqref{g} and~\eqref{4.4} that
  \begin{equation}
|g(z)|\le C_{A'} \exp\left\{-\left|\frac{\theta
x}{4A'}\right|^{1/(1-\beta)} + \left|\frac{\gamma
y}{A}\right|^{1/(1-\beta)}\right\}
 \label{4.5}
  \end{equation}
for every $A'>A$. Since $d(x, U)\le |x|$, we see that the function
$gf$ decreases in the cone $W_1$ if $A<\theta/4$ and the
 inequalities \eqref{4.3} hold in this cone
 with any $B>1+\gamma/A$. On the other hand, we have
\begin{equation}
 d(x, V_1)\ge \frac{\theta}{4}|x| \quad
 \text{for}\quad x\not\in W_1
\label{4.6}
  \end{equation}
by the triangle inequality. Hence $d(x, V_1)\ge \theta d(x, U)/4$
in this region. Since $d(x,U\cup U_1)=\min\{d(x,U), d(x,
U_1\setminus U)\}$ and $U_1\setminus U\subset V_1$, we see that
the inequalities \eqref{4.3}   hold everywhere   if we add the
condition $B\ge 4/\theta$.

Furthermore, the condition $\int g_0(\xi)\d \xi=1$ implies that
$$
(1-g)(z)=\int_{\complement W}g_0(z-\xi)\d \xi.
$$
 Taking
$W_2=\{x\in \oR^d\colon  d(x, V_2)\le \theta|x|/4\}$, we see that
$|x-\xi|\ge \theta |x|/4$ for $x\in  W_2$ and $\xi\in \complement
W$. On the other hand,
\begin{equation}
d(x, V_2)\ge\frac{\theta}{4}|x|\quad \text{for}\quad x\not\in W_2.
\label{4.7}
  \end{equation}
Therefore $(1-g)f\in S^\beta(U\cup U_2)$ provided that
 $A<\theta/4$ as before. This proves the theorem for  $\beta>0$.
\end{proof}
\section{The use of H\"ormander's estimates}

 \begin{proof}[Proof of Theorem 3 for $\beta=0$]
  We first perform a decomposition into smooth functions satisfying
  the bounds  at infinity that are characteristic of
  the elements of $S^0(U\cup U_i)$, and then we restore analyticity.
  Let $W$, $W_1$, $W_2$  be the same auxiliary cones as in the
  previous section. We take an arbitrary nonnegative function
  $\chi_0\in C^\infty_0(\oR^d)$  whose integral is  1 and whose
  support lies in the unit ball, and we set
 $$
 \chi(x)=\int\limits_W\!\chi_0(x-\xi)\,{\rm d}\xi.
 $$
The  argument in Sect.~4 shows that
\begin{gather}
|\chi(x)|e^{d(x,U)}\le Ce^{bd(x,U\cup U_1)}, \label{1}\\
|1-\chi(x)|e^{d(x,U)}\le Ce^{bd(x,U\cup U_2)}, \label{2}\\
\left|\frac{\partial\chi}{\partial x_j}\right|e^{d(x,U)}\le
Ce^{bd(x,U\cup U_1\cup U_2)},\quad j=1,\dots, d,
 \label{3}
\end{gather}
 where we can take  $4/\theta$ for $b$. The situation is
 even simpler than before because $\supp\chi$ is contained
 in the  1-neighborhood  of $W$ and $\chi(x)=0$
at all  points of $W_1$ lying  outside a ball of  sufficiently
large radius  $R$. Inside  the ball, inequality~\eqref{1} holds
with $C=e^R$ and  any $b$. Outside  $W_1$, it holds with $C=1$ and
$b\ge 4/\theta$ by  \eqref{4.6}. Similarly, \eqref{2} follows from
the equation $1-\chi(x)=0$, which holds for all   $x\in W_2$
 outside a sufficiently large ball. The derivatives of
$\chi$ are uniformly bounded, and their supports lie in the
1-neighborhood of the boundary of $W$. At those points of the
supports that lie outside a sufficiently large ball,
both~\eqref{4.6} and  \eqref{4.7} hold, and this yields~\eqref{3}.

We set
 $$
 f=f_1+f_2,\quad f_1(z)=f(z)\chi(x),
 \quad f_2(z)=f(z)(1-\chi(x)).
   $$
 It follows from~\eqref{1} and~\eqref{2} that
   \begin{equation}
 \|f_1\|_{U\cup U_1,b,N}\le C\|f\|_{U, 1, N},\quad
 \|f_2\|_{U\cup U_2,b,N}\le C\|f\|_{U, 1, N},\qquad N=0,1,2,\dots
  \label{4}
  \end{equation}
 To obtain an analytic decomposition, we write
 \begin{equation}
 f= f'_1+ f'_2,\quad f'_1=f_1-\psi,\quad
  f'_2=f_2+\psi
 \notag
\end{equation}
and subject $\psi$ to the equations
\begin{equation}
\frac{\partial\psi}{\partial\bar z_j}=\eta_j ,
 \label{8}
\end{equation}
where
\begin{equation}
\eta_j\stackrel{{\rm def}}{=}f\frac{\partial\chi}{\partial\bar
z_j} =\frac{1}{2}f\frac{\partial\chi}{\partial x_j},\qquad
j=1,\dots,d.
 \label{9}
\end{equation}
By the inequality~\eqref{3*} (which holds even for $B'=B$ in the
case $\beta=0$),    the functions   $\eta_j(z)$ satisfy
\begin{equation}
 |\eta_j(z)|\le C_N \|f\|_{U, 1, N}\,e^{\rho_{{}_{U\cup U_1\cup
 U_2 b,N}}(z)}.
 \label{10}
 \end{equation}
  It remains to show that there is a solution of~\eqref{8}
  with the required behavior at infinity. This can be  done
  using H\"ormander's  $L^2$-estimates. However, the weight functions
  in these estimates are given by exponents of plurisubharmonic functions
while the indicator functions~\eqref{rho} are not of this form.
Therefore we need the following  lemma, which is proved in
Appendix 1.

\begin{lemma} Let $V$ be an open cone in $\oR^d$ and
$B>2edb$. For every function $\eta(z)$, $z\in \oC^n$, satisfying
the inequalities
\begin{equation}
 |\eta(z)|\le C_N\,e^{\rho_{{}_{V, b,N}}(z)},
 \label{11}
 \end{equation}
  there is a plurisubharmonic function $\varrho(z)$ with values in
   $(-\infty,+\infty)$ such that
\begin{equation}
 |\eta(z)|\leq e^{\varrho(z)}\leq
 C'_N\,e^{\rho_{{}_{V, B,N}}(z)}.
 \label{12}
 \end{equation}
\end{lemma}
In our case,  $V=U\cup U_1\cup U_2$ and we apply Lemma 1  to
$\eta=\max |\eta_j|$.  Put
 $$
 \tilde \varrho(z)=2 \varrho(z)+(d+1)\ln(1+|z|^2).
 $$
    Then the functions $\eta_j$  belong to
    $L^2(\oC^d, e^{-\tilde \varrho}{\rm  d}\lambda)$.
   Their definition~\eqref{9} implies that the  compatibility conditions
  $\partial \eta_j/\partial \bar{z}_k =\partial\eta_k/\partial
  \bar{z}_j$
  are fulfilled. By Theorem 15.1.2 of~\cite{H2}, the system
  of equations~\eqref{8} has a solution  $\psi$ such that
   \begin{equation}
   2\int|\psi|^2 e^{-\tilde \varrho}(1+|z|^2)^{-2}{\rm d}\lambda
  \le \int\sum_j|\eta_j|^2 e^{-\tilde \varrho}{\rm d}\lambda.
 \label{13}
   \end{equation}
 It follows from~\eqref{12} and \eqref{13} that
 \begin{equation}
  \psi \in L^2\left(\oC^d,  e^{-2\rho_{V,B,N-d-3}}{\rm
  d}\lambda\right)
    \notag
   \end{equation}
   for any $N$. Combining this with estimates~\eqref{4}, we see that
  $f'_1\in S^{0}(U\cup U_1)$ and   $f'_2\in S^{0}(U\cup U_2)$, as required.
  Theorem 3 is proved.
\end{proof}
 \begin{remark} The  proof of Theorem 3 shows that if $\bar
U_1\cap \bar U_2=\{0\}$, then every element $f\in S^\beta(\{0\})$
may be decomposed as $f=f_1+ f_2$, where $f_i\in S^\beta(U_i)$,
$i=1,2$. This  special case is covered by Theorem 3 with a
slightly changed wording, where the words ``open cones'' are
replaced by ``cones with open projections''. Applying the open
mapping theorem, we see that the topology defined on
$S^\beta(\{0\})$ by the norms $\|\cdot\|_{\{0\},B,N}$, where
$d(x,0)=|x|$, coincides with the inductive topology determined by
the pair of injections
 $S^\beta(U_i)\to S^\beta(\{0\})$, $i=1,2$, as well as with the
 inductive topology determined by the injections
 $S^\beta(U)\to S^\beta(\{0\})$, where $U$ ranges over all
 open cones in $\oR^d$.
\end{remark}

\section{The existence of smallest carrier cones}

 \begin{theorem} For every   continuous linear functional
 on $S^\beta(\oR^d)$, $\beta\ge 0$, there is a unique minimal
 closed carrier cone $K\subset \oR^d$.
 \end{theorem}
 \begin{proof} By Theorem 3, we have
\begin{equation}
  S^\beta(K_1\cap K_2)= S^\beta(K_1)+ S^\beta(K_2),
 \label{14}
   \end{equation}
 for every pair of closed cones in $\oR^d$. Indeed, if
$K_1\cap K_2=\{0\}$, then Remark 1 applies because  there are open
cones $U_i$ such that $K_i\Subset U_i$ and $\bar U_1\cap \bar
U_2=\{0\}$. If $K_1\cap K_2\ne\{0\}$ and $f\in S^\beta(U)$ with
$U\Supset K_1\cap K_2$, then there are $U_i$ satisfying $\Bar U_1
\cap\Bar U_2\Subset U$.

Now we show that~\eqref{14} implies the following dual relation:
\begin{equation}
  S^{\prime\, \beta}(K_1\cap K_2)= S^{\prime\, \beta}(K_1)\cap S^{\prime\,
  \beta}(K_2),
 \label{15}
   \end{equation}
where all spaces are regarded as subspaces of $S^{\prime\,
\beta}(\oR^d)$. The nontrivial part of relation~\eqref{15} states
that if a functional $v\in S^{\prime\, \beta}(\oR^d)$ is carried
by $K_1$ and by $K_2$, then  $K_1\cap K_2$ is also a carrier cone
of $v$. Let $v_i$ be continuous extensions of $v$ to
$S^\beta(K_i)$ and let $f\in S^\beta(K_1\cap K_2)$. Using the
decomposition $f=f_1+f_2$, where $f_i\in S^\beta(K_i)$, we define
an extension of $v$ to $S^\beta(K_1\cap K_2)$ by  $\hat
v(f)=v_1(f_1)+v_2(f_2)$. This extension is well defined. Indeed,
if $f=f'_1+f'_2$ is another decomposition, then
$$
f_1-f'_1=f'_2-f_2\in S^\beta(K_1)\cap S^\beta(K_2)=
S^\beta(K_1\cup K_2)
$$
 and hence $v_1(f_1-f'_1)=v_2(f'_2-f_2)$ because
$S^\beta(\oR^d)$ is dense in $S^\beta(K_1\cup K_2)$ by Theorem~1.
The functional $\hat v$ is obviously continuous in the inductive
topology $\mathcal T$ determined by the injections
$S^\beta(K_i)\to S^\beta(K_1\cap K_2)$, $i=1,2$, and this topology
coincides with the original topology $\tau$ of $S^\beta(K_1\cap
K_2)$ by the open mapping theorem~\cite{G}. Indeed, $\tau$ is not
stronger than $\mathcal T$ and $(S^\beta(K), \mathcal T)$ belongs
to the class $\mathcal{LF}$ because both  spaces $S^\beta(K_i)$
are in this class and $\mathcal T$ coincides with the quotient
topology of the outer sum  $S^\beta(K_1)\oplus S^\beta(K_2)$
modulo a closed subspace (see~\cite{RR}, Ch.~V, Proposition~28).

The relation~\eqref{15} yields an analogous relation for the
intersection of any finite family of closed cones. Then the
existence of a smallest carrier cone for every $v\in
S^{\prime\,\beta}(\oR^d)$ can be established by standard
compactness arguments. Indeed, let $K$ be the intersection of all
carrier cones of $v$ and  let  $U$ be an open cone such that
$U\Supset K$. The projections of the cones complementary to the
carriers cover the compact set $\pr\complement U$, and we can
choose a finite subcovering $\pr\complement K_j$ from this open
 (in the  topology of the unite sphere) covering. Then $\cap_j
K_j\Subset U$. Therefore, the functional $v$ is continuous in the
topology of $S^\beta(U)$, and  $K$ is a carrier cone of $v$. This
proves the theorem.
\end{proof}
Combining~\eqref{14} with the obvious formula
\begin{equation}
  S^\beta(K_1\cup K_2)= S^\beta(K_1)\cap S^\beta(K_2),
 \notag
   \end{equation}
we see that the map $K\to S^\beta(K)$ is a lattice
(anti-)homomorphism from the lattice of closed cones in $\oR^d$ to
the lattice of linear subspaces of $S^\beta(\{0\})$. This is
equivalent to the exactness of the sequence
\begin{equation}
  0  \longrightarrow S^\beta(K_1\cup K_2)
  \stackrel{i}{\longrightarrow}
   S^\beta(K_1) \oplus S^\beta(K_2) \stackrel{s}{\longrightarrow}
   S^\beta(K_1\cap K_2) \longrightarrow 0,
 \label{16}
 \end{equation}
where  $s$ takes each pair of functions $f_{1,2}\in
S^\beta(K_{1,2})$ to the difference of their restrictions to
$K_1\cap  K_2$. As shown above, the sequence~\eqref{16} is even
topologically exact at the term $S^\beta(K_1\cap K_2)$. But we
cannot assert  this for the term $S^\beta(K_1\cup K_2)$. In other
words, we cannot claim that the original topology of this space
coincides with the projective topology determined by the canonical
embeddings into $S^\beta(K_i)$, $i=1,2$. This differs essentially
from the case of the  DFS-spaces $S^\beta_\alpha(K)$ considered
in~\cite{S1}, where the topological exactness of an analogous
sequence  evidently follows from the open mapping theorem, which
applies because any finite sum of DFS spaces and any closed
subspace of a DFS space also belong to this class. However, the
sequence~\eqref{16} is topologically exact in the important event
that $K_1\cap K_2=\{0\}$, because $S^\beta(\{0\})$ is a  DFS
space. Then Theorem 5 of~\cite{S3} shows that every functional
$v\in S^{\prime \,\beta}$ with carrier cone $K_1\cup K_2$, where
$K_1\cap K_2=\{0\}$, can be decomposed into a sum of functionals
$v_i\in S^{\prime\, \beta}(K_i)$, $i=1,2$.

 \section{Nuclearity}

\begin{lemma} For any $B<B'$ and $N>N'$, the natural
   injection   $H^{\beta,B}_N(U)\to H^{\beta,B'}_{N'}(U)$ is a
   Hilbert-Schmidt map.
\end{lemma}
\begin{proof} We  use the fact that holomorphic functions are
pluriharmonic and satisfy the Laplace equation  $\Delta f=0$,
where $\Delta=\sum_j
  (\partial^2/\partial x_j^2+ \partial^2/\partial y_j^2)$.
  As before, we write $L^2(\oC^d, e^{-2\rho}\d\lambda)$
   for the Hilbert space of  complex-valued functions on $\oC^d$
  that are square-integrable with the weight  $\exp\{-2\rho_{U,B,N}\}$,
  where   $\rho_{U,B,N}$ is defined by~\eqref{rho}.
  In what follows we omit the subscripts  $U,B,N$
  and write $\rho'$ for the function specified by $U,B',N'$.
The space $H^{\beta,B}_N(U)$ is a close subspace of $L^2(\oC^d,
e^{-2\rho}d\lambda)$ and hence is separable. We need an auxiliary
function belonging to the space $S^{1-\alpha}_\alpha$, where
$\alpha>\beta$. If $\beta<1/2$, then the function $e^{-t^2}\in
S^{1/2}_{1/2}$ is suitable.  As shown in~\cite{GS}, \S\,\,IV.8,
every space $S_\alpha^{1-\alpha}(\oR)$ with $1/2<\alpha<1$
contains an element of the form $\psi(t^2)$, where $\psi\not\equiv
0$ is an entire function having exponential growth of  order
$1/(2\alpha)$ in the complex plane and   exponential decrease of
the same order along the real semi-axis $t>0$. We assume that
$\psi(0)=1$. Let $p\in \oR^d$, $q\in \oR^d$ and
$\Psi(p,q)=\psi(p^2+q^2)$. According to~\cite{GS}, \S\,\,IV.9, we
have $\Psi\in S_\alpha^{1-\alpha}(\oR^{2d})$ and
\begin{equation}
\Phi(x,y)=\frac{1}{(2\pi)^{2d}}\int e^{-ipx-iqy}\Psi(p,q)\d p\, \d
q \in S^\alpha_{1-\alpha}(\oR^{2d}).
   \notag
      \end{equation}
In particular,  $\Phi$ satisfies the estimate
 $$
 |\Phi(x, y)|\le C
\exp\{-|x/A|^{1/(1-\alpha)}-|y/A|^{1/(1-\alpha)}\}
$$
with some $A>0$. Therefore  the convolution $\Phi*f$ exists for
any function $f\in L^2(\oC^d, e^{-2\rho}\d\lambda)$ if
$\alpha>\beta$. Let $\rho_1=\rho_{U,B_1,N}$, where $B_1>B$.
Applying the triangle inequality to each term of  $\rho$, we find
that
\begin{equation}
\int |\Phi(x',y')|^2 e^{2\rho(x-x',y-y')}\d \lambda'\le
C'e^{2\rho_1(x,y)}.
   \label{7.1}
      \end{equation}
Using next the Cauchy-Schwarz-Bunyakovskii inequality, we obtain
\begin{equation}
 |(\Phi*f)(x,y)|\le C'\|f\|_{U,B,N} e^{\rho_1(x,y)}.
   \notag
      \end{equation}
Choosing  $B_1<B'$, we see that the correspondence $f\to \Phi*f$
is a continuous map from $L^2(\oC^d, e^{-2\rho}\d\lambda)$ to
$L^2(\oC^d, e^{-2\rho'}\d\lambda)$. Moreover, it is a
Hilbert-Schmidt map. Indeed, the multiplication by $e^{-\rho}$ is
an isometry from $L^2(\oC^d, e^{-2\rho}\d\lambda)$ onto
$L^2(\oC^d)$ and the map in question belongs to the
Hilbert-Schmidt class if and only if the integral operator on
$L^2(\oC^d)$ with  kernel
$e^{-\rho'(x',y')}\Phi(x-x',y-y')e^{\rho(x,y)}$ is in the same
class, that is, if the kernel is square-integrable, and this  is
ensured by the estimate~\eqref{7.1}.

On the other hand, the map $f\to \Phi*f$ is identified with the
infinite-order differential operator
$\psi(-\Delta)=1+\sum_{k\ge1}c_k\Delta^k$ if $f$ is treated as a
generalized function defined on appropriate test functions.
According to~\cite{GS}, such an operator is well defined on any
space $S^{\beta'}_{1-\beta'}(\oR^{2d})$, where $\beta'<\alpha$. If
$\beta'>\beta$, then all elements of $L^2(\oC^d,
e^{-2\rho}\d\lambda)$ are integrable with  test functions in this
space. With this choice of $\beta'$, we have the chain of
identities
 $$
 (f, \psi(-\Delta)\varphi)=\lim_{n\to\infty}\left(f,
 \left(1+\sum^n_{k=1}c_k\Delta^k\right)\varphi\right)=(2\pi)^{-2d}(\tilde f,
 \Psi\tilde\varphi)= (f,\Phi*\varphi),
 $$
where $\varphi$ is any element of
$S^{\beta'}_{1-\beta'}(\oR^{2d})$. In particular,
$(\Phi*f,\varphi)=(f,\varphi)$ for all $f\in H^{\beta,B}_N(U)$. It
follows that $\Phi*f=f$ because $S^{\beta'}_{1-\beta'}(\oR^{2d})$
has a sufficiently large stock of functions (see~\cite{GS}). This
proves Lemma~2.
\end{proof}

  \begin{theorem}   The spaces $S^{\beta,b}(U)$ and $S^\beta(U)$
   are nuclear for any open cone $U\subset \oR^d$.
   The spaces  $S^\beta(K)$ associated with closed cones are also nuclear.
\end{theorem}
 \begin{proof} The statement for $S^{\beta,b}(U)$ follows
  immediately from Lemma 2 because the composite of two
  Hilbert-Schmidt maps is  nuclear  and the projective
  limit of a sequence of Hilbert spaces with nuclear connecting
  maps is a nuclear Fr\'echet space. The statement about $S^\beta(U)$
  and $S^\beta(K)$ follows from the heredity properties
  of  inductive limits of countable families of nuclear spaces
 (see~\cite{Sch}).
\end{proof}
\begin{corollary} The spaces $S^{\beta,b}(U)$ and $S^\beta(U)$
 are reflexive. Moreover, they are Montel spaces.
\end{corollary}
Indeed, they are complete and barrelled, and every nuclear space
with these properties is a Montel space (see~\cite{Sch}, Ch.~IV,
Exercise~19). It is still an open question  whether the spaces
$S^\beta(K)$ over closed cones have these properties. But their
completions certainly have them.

\section{Kernel theorems}

    If $E_1$ and  $E_2$ are locally convex spaces (LCS), then their
    (algebraic)  tensor product equipped with the projective
    topology $\tau_\pi$ is denoted by  $E_1\otimes_\pi E_2$,
    and the same product  with the
    inductive topology $\tau_\iota$ is denoted by
    $E_1\otimes_\iota E_2$. If   $E_1$ and $E_2$ are Hilbert spaces,
    then we write  $E_1\otimes_{\rm H} E_2$ for their tensor product
    equipped with the natural inner product. The completion of each
    of these spaces is denoted by a ``hat''  over the tensor product symbol.

  \begin{lemma} Let $U_i$ be open cones in $\oR^{d_i}$, $i=1,2$.
  Then there is a canonical isomorphism
    \begin{equation}
   H^{\beta,B}_N (U_1)\mathbin{\hat{\otimes}_{\rm H}}
    H^{\beta,B}_N (U_2) \simeq H^{\beta,B}_N (U_1\times U_2),
  \notag
  \end{equation}
defined by  identifying $f_1\otimes f_2$ with the function
$f_1(z_1) f_2(z_2)$.
\end{lemma}

 \begin{proof} We use the same line of reasoning
as in the case of square-integrable functions. (This case is
considered, for example,  in~\cite{He}.) By the
property~\eqref{m}, all function of the form $f_1(z_1)f_2(z_2)$,
where $f_1\in H^{\beta,B}_N (U_1)$ and $f_2\in H^{\beta,B}_N
(U_2)$, belong to the space $H^{\beta,B}_N (U_1\times U_2)$, and
their linear span is identified with $H^{\beta,B}_N (U_1)\otimes
H^{\beta,B}_N (U_2)$. The natural inner product on a tensor
product of   Hilbert spaces is defined by
   $$
   \langle f_1\otimes f_2, g_1\otimes g_2\rangle=\langle f_1, g_1\rangle
   \langle  f_2,  g_2\rangle
   $$
 with  subsequent extension by linearity. In our case
 it obviously coincides with the inner product induced by  that
  of $H^{\beta,B}_N (U_1\times U_2)$. If $\{f_j\}$
and $\{g_k\}$ are bases in the spaces whose tensor product is
being formed, then $\{f_j(z_1)g_k(z_2)\}$ is an orthonormal system
in $H^{\beta,B}_N (U_1\times U_2)$ and  Fubini's theorem
 immediately shows that this system is total. Therefore the
completion of the tensor product coincides with $H^{\beta,B}_N
(U_1\times U_2)$.
\end{proof}
 \begin{lemma} Let  $h_1\colon E_1\to F_1$ and
$h_2\colon E_2\to F_2$  be  Hilbert-Schmidt maps between Hilbert
spaces.  Then the map
   \begin{equation}
   E_1\otimes_{\rm H} E_2\stackrel{h_1\otimes h_2}{\longrightarrow}
   F_1\otimes_\pi F_2
  \notag
   \end{equation}
   is continuous.
   \end{lemma}
\begin{proof} We assume that all the spaces are separable because
this is the case in the applications below, although this lemma
holds in the general case as well. A map $h\colon E\to F$ belongs
to the Hilbert-Schmidt class if and only if
$$
\|h\|_2   \stackrel{\mathrm{def}}{=}\left(\sum_j
   \|he_j\|^2 \right)^{1/2}<+\infty
$$
    for some (and thus for every)
   orthonormal basis $\{e_j\}$ in $E$. According to~\cite{Sch},
    Ch.~III, \S\,\,6.3, the projective topology on $F_1\otimes F_2$
    is determined by  the tensor product  of the norms on
$F_1$ and  $F_2$. Denoting this product by $\|\cdot\|_\pi$, we
recall that it is a cross-norm, that is, $\|f_1\otimes
f_2\|_\pi=\|f_1\|\cdot\|f_2\|$.  Moreover, it is stronger than any
other cross-norm. In particular, it is stronger   than the Hilbert
norm determined by the  inner product. Let  $\{e_j^1\}$
    and $\{e_k^2\}$ be orthonormal bases in $E_1$ and
   $E_2$ respectively. Then  $\{e_j^1\otimes e_k^2\}$
   is an orthonormal basis in
    $E_1\hat{\otimes}_{\rm H}E_2$
   and  every  element  $g$ of this space can be written as
   $g=\sum\lambda_{jk}\,   e_j^1\otimes e_k^2$.
    Using the cross-property of  $\|\cdot\|_\pi$, the
    Cauchy-Schwarz-Bunyakovskii  inequality and  Parseval's identity
 $\|g\|^2=\sum|\lambda_{jk}|^2$, we obtain
   \begin{equation}
   \|(h_1\otimes h_2)\left(\sum_{jk\le n}\lambda_{jk}\,
     e_j^1\otimes e_k^2\right)\|_\pi \leq \sum_{j,k\le n}|\lambda_{jk}|\,
   \|h_1(e_j^1)\|\, \|h_2(e_k^2)\|\leq
   \|g\|\,\|h_1\|_2\,\|h_2\|_2.
    \notag
   \end{equation}
It follows that the family $\lambda_{jk}\, h_1(e_j^1)\otimes
h_2(e_k^2)$ of elements of the Banach space $F_1\hat{\otimes}_\pi
F_2$ is absolutely summable. Hence we have defined a continuous
map $E_1\hat{\otimes}_{\rm H} E_2\to F_1\hat{\otimes}_\pi F_2$.
This map coincides with $h_1\otimes h_2$ on the basis elements and
hence on all elements of $E_1\otimes E_2$ because the canonical
bilinear map of $F_1\times F_2$ to $F_1\hat{\otimes}_{\rm H}F_2$
is continuous. The lemma is  proved.
\end{proof}
 \begin{lemma}  The space $S^{\beta,b}(U_1)\otimes
  S^{\beta,b}(U_2)$ is dense in
$S^{\beta,b}(U_1\times U_2)$ for every pair of open cones
$U_i\subset\oR^{d_i}$, $i=1,2$, and  any $0\le\beta<1$, $b>0$.
\end{lemma}
\begin{proof} Let $1<\alpha<2$ and let $g$ be an entire function
on $\oC^{d_1+d_2}$ satisfying
\begin{equation}
|g(z)|\le C\exp\left\{-\sum_j|x_j|^{1/\alpha} +
|by|^{1/(1-\beta)}\right\}
   \notag
   \end{equation}
  and such that $g(0)=1$. We take $f\in S^{\beta,b}(U_1\times U_2)$
  and consider the sequence
     $f_\nu(z)=f((1-1/\nu)z)g(z/(2\nu))$, $\nu=1,2,\dots$. Setting
     $\epsilon=1/\nu$, $p=1/(1-\beta)$, and using the inequalities
  $1\ge (1-\epsilon)^p+\epsilon^p>(1-\epsilon)^p+(\epsilon/2)^p$,
     we  easily verify that $f_\nu$ is bounded in  each of
the norms of $S^{\beta,b}(U_1\times U_2)$. Therefore $f_\nu\to f$
in the topology of this space because it is a Montel space and its
topology is stronger than the  topology of  pointwise convergence.
Let $H_{1/2}(U)$ denote the Hilbert space of entire functions
belonging to $L^2(\oC^d, e^{-2\rho_{1/2}}\d\lambda)$, where
$$
\rho_{1/2}=\exp\left\{-\sum_j|x_j|^{1/2}+  d(bx,U)^{1/(1-\beta)} +
|by|^{1/(1-\beta)}\right\}.
$$
Clearly, it is continuously embedded in $S^{\beta,b}(U)$. All the
functions $f_\nu$ are contained in $H_{1/2}(U_1\times U_2)$.
  Indeed, if $B>b$ and is sufficiently close to $b$, then
 \begin{equation}
  \|f_\nu\|_{1/2}\le C_\nu\|f\|_{U_1\times U_2,B,0}.
   \notag
   \end{equation}
The  arguments used in the proof of Lemma 3 show that
$H_{1/2}(U_1)\otimes H_{1/2}(U_2)$ is dense in $H_{1/2}(U_1\times
U_2)$. Therefore every function $f_\nu$ can be approximated by
elements of the tensor product in a metric  stronger than that of
$S^{\beta,b}(U_1\times U_2)$. This proves the lemma.
\end{proof}
\begin{theorem}
Let $U_1$, $U_2$ be open cones in $\oR^{d_1}$, $\oR^{d_2}$.  Then
there are the canonical isomorphisms
\begin{gather}
        S^{\beta, b} (U_1)\mathbin{\hat{\otimes}_\iota} S^{\beta, b} (U_2)
     \simeq S^{\beta, b} (U_1\times U_2)
 \label{8.1}\\
   S^\beta (U_1)\mathbin{\hat{\otimes}_\iota} S^\beta (U_2) \simeq
S^\beta (U_1\times U_2).
  \label{8.2}
  \end{gather}
\end{theorem}
\begin{proof} The topologies  $\tau_\iota$ and $\tau_\pi$
coincide
  on tensor products of  Fr\'echet spaces (see~\cite{Sch}, Ch.~III,
   \S\,\,6.5), and  Lemmas 2--4 show that the topology $\tau_\pi$ on
    $S^{\beta, b} (U_1)\otimes S^{\beta, b} (U_2)$ coincides with the topology
 induced by that of $S^{\beta, b} (U_1\times U_2)$ because the systems of
 norms determining these topologies are equivalent. By Lemma 5,
  the natural injection
    $ S^{\beta, b} (U_1)\mathbin{\otimes_\iota} S^{\beta, b} (U_2)
     \to S^{\beta, b} (U_1\times U_2)$ has a unique extension
     to the completion of the tensor product and this extension is
an isomorphism. This proves~\eqref{8.1}. Isomorphism~\eqref{8.2}
follows     from~\eqref{8.1} because of  two facts.    First, if
$E_\nu$ and $F_\nu$ are injective sequences of locally convex
spaces and their inductive limits are
 Hausdorff spaces\footnote{The very definition of a LCS usually requires
 the space to be  Hausdorff. But  this property can be lost
 after taking an inductive limit.} then
\begin{equation}
(\varinjlim E_\nu)\otimes_\iota (\varinjlim F_\nu)=\varinjlim
(E_\nu\otimes_\iota F_\nu).
    \label{8.3}
   \end{equation}
   Second, if $G_\nu$ is an injective sequence
    of locally convex spaces and the limit $\varinjlim \hat G_\nu$ of their
    completions is a Hausdorff space, then
\begin{equation}
\widehat{\varinjlim  G_\nu}=\widehat{\varinjlim \hat G_\nu}.
    \label{8.4}
     \end{equation}
We set $E_\nu=S^{\beta,\nu}(U_1)$, $F_\nu=S^{\beta,\nu}(U_2)$,
$G_\nu=E_\nu\otimes_\iota F_\nu$, and successively
use~\eqref{8.3}, \eqref{8.4}, \eqref{8.1}. Since
$S^\beta(U_1\times U_2)$ is  complete, we get~\eqref{8.2}.

The relation~\eqref{8.3} is actually a part of Proposition~14 in
Ch.~I of~\cite{G}. We note that the topology $\tau_\iota$ on a
tensor product of locally convex spaces is certainly Hausdorff
because it is stronger than the topology  $\tau_\pi$, which is
 Hausdorff by~\cite{RR},  Ch.~VII, Proposition~8. The proof
of~\eqref{8.3} consists of using the definition of
 $\tau_\iota$ as the topology of uniform convergence on
 separately equicontinuous sets of bilinear forms and   noting  that a set
  of bilinear forms on $(\varinjlim E_\nu)\times (\varinjlim F_\nu)$ is
separately equicontinuous if and only if the same is true for the
sets of their restrictions to each of $E_\nu\times F_\nu$.
  To prove~\eqref{8.4}, we start by  noting that the continuous
  injections $u_{\mu\nu}\colon G_\mu\to G_\nu$, $\nu>\mu$,
generate continuous maps $\hat u_{\mu\nu}\colon \hat G_\mu\to \hat
G_\nu$ which still satisfy the chain rule $\hat u_{\nu\mu}\circ
\hat u_{\mu\lambda}=\hat u_{\nu\lambda}$. Therefore the space
$\varinjlim \hat G_\nu$ is well defined. For every $\nu$, we have
the continuous map $G_\nu\to\varinjlim \hat G_\nu$. It is
injective because  the restriction of $\hat u_{\mu\nu}$ to $G_\nu$
is one-to-one whenever $\mu>\nu$. These injections determine a
continuous injection
\begin{equation}
\varinjlim G_\nu\to \varinjlim \hat G_\nu.
 \label{8.5}
     \end{equation}
In particular, if $\varinjlim \hat G_\nu$ is  Hausdorff, then so
is $\varinjlim G_\nu$. On the other hand, the injections $G_\nu\to
\varinjlim G_\nu$ extend to maps $\hat G_\nu
\to\widehat{\varinjlim  G_\nu}$ and generate a continuous map
\begin{equation}
\varinjlim \hat G_\nu\to \widehat{\varinjlim  G_\nu}.
 \notag
     \end{equation}
The composite of this map and~\eqref{8.5} is the canonical
embedding of the space $\varinjlim G_\nu$ into its completion.
Therefore the topology of  $\varinjlim G_\nu$ coincides with the
topology induced by that of $\varinjlim \hat G_\nu$. Since the
image of $\varinjlim G_\nu$ is dense in $\varinjlim \hat G_\nu$,
the injection~\eqref{8.5} extends to an isomorphism, which
completes the proof of Theorem 6.
\end{proof}

To every closed cone in  $\oR^{d_1+d_2}$ having the product
structure $K_1\times K_2$ with $K_i\subset\oR^{d_i}$ we assign the
space
\begin{equation}
S^\beta(K_1, K_2)=\varinjlim_{U_1,U_2} S^\beta(U_1\times U_2),
 \label{8.6}
     \end{equation}
where the $U_i$ are open cones in $\oR^{d_i}$ such that
$U_i\Supset K_i$, $i=1,2$.

 The existence of canonical embeddings $S^\beta(U_1\times U_2)\to
S^\beta(\{0\})$ enables us to interpret~\eqref{8.6} as an
inductive limit in the sense of the definition  given
in~\cite{G,RR}. However,  if one (or both) of  the cones $K_i$ is
degenerate, then the set of open cones involved in the limit is
not directed. Directed sets are sometimes more convenient to use,
and the situation
  can  easily be remedied by taking a limit over all cones
  with open projections in $\oR^{d_i}$ instead of open cones.
  In other words, we can simply add the cone $\{0\}$ to the set
  of open cones  in each  $\oR^{d_i}$. This agrees with
  what was said in Sect.~2. Spaces associated with
    relatively open cones  $U_1\times\{0\}$ and $\{0\}\times
U_2$ are defined in the same manner as those associated with open
cones, and Theorem~6 (as well as Theorem~1) can be immediately
extended to them. In particular, if $\{0\}$ is a degenerate cone
in $\oR^{d_1}$, then
$$
S^\beta (\{0\})\mathbin{\hat{\otimes}_\iota} S^\beta (U_2) \simeq
S^\beta (\{0\}\times U_2).
$$ Such a
unification was used in~\cite{Sm}, where an analogue
of~\eqref{8.6}  was proposed for the DFS-spaces $S^\beta_\alpha$.
We emphasize that the addition of $\{0\}$ to the set of open cones
leaves inductive limit~\eqref{8.6} unchanged. For instance, let
$K_1=\{0\}$. Choose any two open cones in $\oR^{d_1}$ with
disjoint closures of their projections (say, the positive and
negative orthants $U_+$ and $U_-$). Theorem 3 shows that every
element of $S^\beta(\{0\}\times U_2)$ can be written as a sum of
functions belonging to $S^\beta(U_+\times V_2)$ and
$S^\beta(U_-\times V_2)$, where $U_2\Supset V_2\Supset K_2$. Using
the open mapping theorem, we see that the inductive topology
 on $S^\beta(\{0\}, K_2)$ with respect to the family of subspaces
 $S^\beta(U_1\times U_2)$, where $U_i\Supset K_i$, coincides with
that determined by the subfamily $S^\beta(U_{\pm}\times U_2)$,
where $U_2\Supset K_2$. It also coincides with the inductive limit
topology with respect to the increasing family
$S^\beta(\{0\}\times U_2)$.

 \begin{theorem} Let $K_i$ be a closed cone in
 $\oR^{d_i}$, $i=1,2$. Every separately continuous bilinear form
  $w$ on $S^\beta(K_1)\times S^\beta(K_2)$ is uniquely representable as
   \begin{equation}
   w(f,g) = (v,f\otimes g),
   \notag
   \end{equation}
   where $v$ is a continuous linear functional on   $S^\beta(K_1, K_2)$.
  \end{theorem}

\begin{proof}   By the isomorphism~\eqref{8.2}, the restriction of
$w$ to $S^\beta(U_1)\times S^\beta(U_2)$, where $K_i\Subset U_i$,
uniquely determines a continuous linear functional on
$S^\beta(U_1\times U_2)$. If both the cones $K_1$, $K_2$ are
nondegenerate, then the family of neighborhoods  $U_1\times U_2$
is decreasing and hence we have defined a linear functional on
$S^\beta(K_1, K_2)$, which is continuous by the definition of the
inductive topology. The same argument works in the degenerate case
if we use the above-stated unification of definition~\eqref{8.6}
and the corresponding generalization of Theorem 6. Another way is
to use Theorem 3.  For instance, let $K_1=\{0\}$ as above. If
$v_{\pm}$ are the functionals generated by $w$ on
$S^\beta(U_{\pm}\times V_2)$, then we can use the decomposition
$f=f_+ +f_-$ with $f_{\pm}\in S^\beta(U_{\pm}\times V_2)$ and
define $v$ by $v(f)=v_+(f_+)+v_-(f_-)$. This functional is well
defined because $v_+$ and $v_-$ coincide on $S^\beta(U_+\times
V_2)\cap S^\beta(U_-\times V_2)=S^\beta((U_+\cup U_-)\times V_2)$
by~\eqref{8.2}.
\end{proof}
\begin{corollary} The space of all separately continuous bilinear
forms on $S^\beta(K_1)\times S^\beta(K_2)$ can be identified with
the space $S^{\prime\,\beta}(K_1, K_2)$, which is the dual of
$S^\beta(K_1, K_2)$.
\end{corollary}

\section{Carrier cones of multilinear forms}

Now we  turn  to the carrier cones of  functionals  generated by
multilinear forms on $S^\beta(\oR^{d_1})\times\dots \times
S^\beta(\oR^{d_n})$ with given carrier cones for every argument.
Our main concern is about the extension of such forms to larger
spaces. The general theory of extension of multilinear forms is
based on the notion of hypocontinuity~\cite{Sch}, but we need only
the following simple lemma.

\begin{lemma}
Let $L$ be a sequentially dense subspace of a locally convex space
$E_1$, and let $E_2$ be a barrelled space. Then every   separately
continuous form bilinear $w$  on $L\times E_2$ has a unique
extension to $E_1\times E_2$ which is bilinear and separately
continuous.
\end{lemma}

 \begin{proof} For every fixed  $g\in E_2$, the form
$w(f,g)$  extends uniquely to $E_1$ by continuity. We  denote this
extension by $\hat w$. We  must verify that this functional is
linear and continuous in $g$ for every fixed $f\in E_1$. Choose a
sequence $f_\nu\in L$ that converges to $f$ and denote the
corresponding elements of $E'_2$ by $\mathbf f_\nu$. Then $\mathbf
f_\nu(g)=w(f_\nu, g)$. It is well known (see~\cite{RR}
or~\cite{Sch}) that if $E_2$ is barrelled, then the pointwise
convergence of the sequence $\mathbf f_\nu\in E'_2$ implies that
its limit $\mathbf f$  belongs to $E'_2$, i.e., is linear and
continuous. We have $\mathbf f(g)=\hat w(f,g)$. This proves the
lemma.
\end{proof}
\begin{definition} Let    $w$ be a multilinear separately
continuous form on
 $S^\beta(\oR^{d_1})\times\dots \times S^\beta(\oR^{d_n})$,  and let
 $K=K_1\times\dots\times K_n$, where $K_j$ is a closed cones in
$\oR^{d_j}$, $j=1,\dots, n$. We say that $K$ is a {\it  carrier
cone} of $w$, if every $K_j$ is a carrier cone for all linear
functionals defined on $S^\beta(\oR^{d_j})$ by $w(f_1,\dots,f_n)$
with fixed  $f_i\in S^0(\oR^{d_i})$,    $i\neq j$.
\end{definition}
\begin{theorem} Let $K=K_1\times\dots\times K_n$, where
$K_j\subset\oR^{d_j}$. If $K$ is a carrier cone of a multilinear
separately continuous form $w$  on $S^\beta(\oR^{d_1})\times\dots
\times S^\beta(\oR^{d_n})$, then $K$ is also a carrier cone of the
functional $v\in S^{\prime \beta}(\oR^{d_1+\dots+d_n})$ generated
by this form.
\end{theorem}

\begin{proof} For simplicity we assume that $d_j=d$ for all $j$.
It suffices to show that every cone
 $\oR^{d(j-1)}\times K_j\times\oR^{d(n-j)}$ is a carrier of $v$
 (so that their intersection is also a carrier by Theorem 4).
 We put $j=1$ without loss of generality.
Suppose that $n=2$  and  $U_1$ is an open cone in $\oR^{d}$ such
that $K_1\Subset U_1$. By Theorem 1, each element  $f\in
S^\beta(U_1)$ can be approximated by elements of
$S^\beta(\oR^{d})$ in the metric of some $S^{\beta,b}(U_1)$ (where
$b$ is dependent on $f$).
 Applying Lemma 6, we extend  $w$  to a
bilinear separately
  continuous form on $S^\beta(U_1)\times S^\beta(\oR^{d})$. By Theorem  6,
this form in turn determines a linear continuous functional on
$S^\beta(U_1\times\oR^{d})$, which is an extension of  $v$ because
$S^\beta(\oR^{d})\otimes S^\beta(\oR^{d})$ is dense in
$S^\beta(\oR^{2d})$. This proves Theorem 8  for bilinear forms
because the intersection of all cones $\bar U_1\times \oR^d$ is
equal to $K_1\times \oR^d$.

Now we use induction on $n$. We regard the $n$-linear form $w$
($n>2$)  as a bilinear form on $L\times E_2$, where
$L=\stackrel{n-1}{\otimes}\!S^\beta(\oR^d)$ and
$E_2=S^\beta(\oR^{d})$.
  By the inductive hypothesis, this form is separately
  continuous  if $L$ is given the topology induced by that of
   $E_1=S^\beta(U_1\times\oR^{d(n-1)})$. The subspace $L$ is sequentially
   dense in $E_1$ because every element of
 $S^\beta(U_1\times\oR^{d(n-1)})$ can be approximated  in the metric of
 $S^{\beta,b}(U_1\times\oR^{d(n-1)})$ by elements of some
  $S^{\beta,b}(\oR^{dn})$, which can in turn be
 approximated by elements of $\stackrel{n-1}{\otimes}\!S^{\beta,b}(\oR^d)$
  in the (stronger) metric of  $S^{\beta,b}(\oR^{dn})$. Again
 using Lemma~6, we conclude that $v$  has a continuous extension to
 $S^\beta(U_1\times\oR^{d(n-1)})$. This  proves the theorem.
\end{proof}
Theorems 7 and 8 raise the question of the relation between the
spaces $S^\beta(K_1\times K_2)$ and $S^\beta(K_1, K_2)$. Clearly,
we have
 $$
 S^\beta(K_1\times K_2)\subset S^\beta(K_1, K_2)
 $$
because $U\Supset K_1\times K_2$ implies that $U\supset U_1\times
U_2$, where $U_i=U\cap\oR^{d_i}$ if $K_i\ne\{0\}$ and $U_i=\{0\}$
otherwise. As a rule, this inclusion is  strict.

  \begin{theorem} The spaces $S^\beta(K_1\times K_2)$
and
  $S^\beta(K_1, K_2)$  coincide only  if $K_1= \oR^{d_1}$ and
  $K_2=\oR^{d_2}$ or if  both these cones are degenerate.
  In all other cases, these spaces are distinct  and have
  different dual spaces.
  \end{theorem}

A proof of Theorem 9 is given  in Appendix~2. It gives an idea of
the stock of functions in the space
  $S^\beta$ that ensures the angular localizability of the
  functionals belonging to its dual.
An obvious generalization of definition~\eqref{8.6} is
  \begin{equation}
 S^\beta(K_1,\dots, K_n)=\varinjlim_{U_1,\dots U_n}S^\beta(U_1\times\dots\times
 U_n),
  \label{9.1}
  \end{equation}
  where  $K_j\Subset U_j$.
   If a functional $v\in S^{\prime\, \beta}(\oR^{d_1+\dots+d_n})$
 admits a continuous extension to the space~\eqref{9.1}, then we say that
 it is {\it strongly carried} by the cone
  $K_1\times\dots\times K_n$. These are precisely those functionals
that  are generated by multilinear forms carried by
$K_1\times\dots\times K_n$.
  To prove this, we need another decomposition theorem.

 \begin{theorem} Let $K_j$ be a closed cone in
  $\oR^{d_j}$, $j=1,\dots, n$.  Every
  $f\in S^\beta(K_1,\dots, K_n)$ can be  decomposed as
 $f=f_1+\dots+ f_n$, where
  $f_j\in S^\beta(\oR^{d_1+\dots +d_{j-1}},
 K_j,\oR^{d_{j+1}\dots +d_n})$.
 \end{theorem}

\begin{proof} This is basically the same as that of Lemma~3
in~\cite{Sm} on the DFS  spaces $S^\beta_\alpha(K_1,\dots, K_n)$.
Suppose that $n=2$. If $K_1$ and  $K_2$ are nondegenerate, then
there are open cones $V_1$ and $V_2$ such that $f\in
S^\beta(V_1\times V_2)$. In the degenerate case,  $f$  is a sum of
elements of such spaces. It suffices to consider the first case.
We choose open cones $V'_j$ such that
 $V_j\Supset V'_j\Supset K_j$ and use the notation $U=V_1\times V_2$,
 $U_1=V'_1\times \complement \bar V_2$,
 $U_2=\complement \bar V_1\times V'_2$.  Then $\bar U_1\cap \bar U_2= \{0\}$
 and it follows from Theorem 3 that $f=f_1+ f_2$, where $f_{1,2}\in
 S^\beta(U\cup U_{1,2})$.  Furthermore,
 $\overline{U\cup U_1} =(\bar V_1\times \bar V_2)\cup\bar
 U_1\supset V'_1\times \oR^{d_2}$. Since the inclusion
 $\overline{W}\supset W'$ implies that $S^\beta(W)\subset S^\beta(W')$,
we conclude that
    $f_1\in S^\beta(V'_1\times \oR^{d_2})\subset S^\beta(K_1,\oR^{d_2})$.
Similarly, $f_2\in  S^\beta(\oR^{d_1},K_2)$.  The same argument
shows that every element of
  $S^\beta(V_1\times\dots\times V_n)$ (with $m$ cones $V_j \Supset K_j$
being different from $\oR^{d_j}$) is representable as a sum of two
functions belonging to the spaces $S^\beta(V'_1\times\dots\times
V'_n)$, where
 $V_j\Supset V'_j\Supset K_j$ and $m-1$ cones
 $V'_j$  are different from $\oR^{d_j}$. Indeed, let the cones
 $V_j\ne\oR^{d_j}$ occupy the first place. For the rest, we set
 $V'_j=V_j=\oR^{d_j}$ and now use the notation
   $U=V_1\times\dots\times V_n$, $U_1=V'_1\times \complement \bar V_2\times
 V'_3\times\dots\times V'_n$, $U_2=\complement \bar V_1\times V'_2\times
 V'_3\times\dots\times V'_n$.  Then $\bar U_1\cap \bar U_2\Subset
 U$ and Theorem~3 again  applies.
\end{proof}

 \begin{theorem} Let $K_j$ be a closed cone in
   $\oR^{d_j}$, $j=1,\dots n$. A functional $v\in S^{\prime\,
   \beta}(\oR^{d_1+\dots+d_n})$ is strongly carried by
 the cone $K_1\times\dots\times K_n$ if and only if $v$ is
 generated by a multilinear separately  continuous form
 on $S^\beta(\oR^{d_1})\times\dots \times S^\beta(\oR^{d_n})$
  that is carried by this cone.
  \end{theorem}

\begin{proof} Let  $v\in S^{\prime\,\beta}(K_1,\dots, K_n)$,
$U_j\Supset K_j$ and $f_j\in S^\beta(U_j)$. It is clear
from~\eqref{m} that the map
$$
f_j\to f_1\otimes\dots\otimes f_j \otimes\dots\otimes f_n
 $$
  from
 $S^\beta(U_j)$ to  $S^\beta(\oR^{d(j-1)}\times U_j\times\oR^{d(n-j)})$
  is continuous. Therefore, the multilinear form corresponding to $v$ is
  certainly carried by $K_1\times\dots\times K_n$.
 To prove the converse, we again set $d_j=d$  for simplicity.
 If some of the  $K_j$ are degenerate, then it is convenient to use
  the unification of
  definitions~\eqref{8.6}, \eqref{9.1} mentioned in Sect.~8.
   In the proof of Theorem~8 we saw that
 $v\in S^{\prime \beta}(\oR^{d(j-1)}, K_j,\oR^{d(n-j)})$ for any
 $j$. Let us show that there is a continuous extension
 $\hat v$ to the space
$E=S^\beta(K_1,K_2,\oR^{d(n-2)})$. By Theorem~10, this space is
the sum of  the two subspaces $L_1=S^\beta(K_1, \oR^{d(n-1)})$ and
$L_2=S^\beta(\oR^{d}, K_2,\oR^{d(n-2)})$. If $f=f_1+f_2$, we set
 $\hat v(f)=v_1(f_1)+v_2(f_2)$. This extension is well defined because
 $v_1$ and $v_2$ coincide on $L_1\cap L_2$. Indeed, this
 intersection is the inductive limit of the increasing family of
 spaces
 $S^\beta(U)$, where $U$ is the union of cones
 $U_1\times\oR^{d(n-1)}$ and $\oR^{d}\times U_2\times\oR^{d(n-2)}$
 with $U_1\Supset K_1$ and $U_2\Supset K_2$. By Theorem~1,  the space
$S^\beta(\oR^{dn})$ is dense in this intersection, whose topology
is stronger than the topologies of  $L_1$ and $L_2$. The
functional $\hat v$ is obviously continuous in the inductive
topology determined by the injections $L_i\to E$, $i=1,2$, and
this topology coincides with the original topology of  $E$ by the
open mapping theorem. Applying the same arguments to the triple
$E=S^\beta(K_1,K_2, K_3,\oR^{d(n-3)})$,
$L_1=S^\beta(K_1,K_2,\oR^{d(n-2)})$, $L_2=S^\beta(\oR^{2d},
K_3,\oR^{d(n-3)})$ and so on, we complete the proof after finitely
many steps.
\end{proof}

\section{A Paley-Wiener-Schwartz-type theorem}

 Let $V$ be an open cone in $\oR^d$ and let
 $V^*=\{x\colon  x\eta\ge 0, \forall  \eta\in V\}$ be its dual  cone.
 As shown in~\cite{FS},  the Laplace transformation maps the space
  $S^{\prime \beta}(V^*)$, $\beta>0$, onto an algebra of analytic functions
  defined on the tubular domain $T^V=\oR^{dn}+iV$ and satisfying
  certain bounds on their growth
  near the real boundary of the domain of analyticity and at  infinity.  An
  analogous theorem was proved in~\cite{S} for the  class
   $S^{\prime 0}$, which requires more
  sophisticated reasoning. Theorem~11 enables us to extend
  these results to the spaces $S^{\prime\,\beta}(V_1^*,\dots, V_n^*)$.

Let $\beta>0$, let $V_j$ be open cones in $\oR^{d_j}$, $j=1,\dots,
n$, and let  $V=V_1\times\dots\times V_n$. We denote by $\mathcal
A_\beta(V_1,\dots, V_n)$ the space of functions analytic in the
domain $T^V=\oR^{dn}+iV$ and satisfying the condition
 \begin{equation}
  |\mathbf u(\zeta)|\leq
 C_{\epsilon,W_1,\dots,W_n}\prod_{j=1}^n |\I \zeta_j|^{-N}
 \exp\{\epsilon |\R \zeta_j|^{1/\beta}\},\qquad
   \I \zeta_j\in W_j,\,  j=1,\dots,n
 \label{10.1}
 \end{equation}
for any $\epsilon>0$,   any cones $W_j\Subset V_j$ and  some $N$
depending on  $\epsilon$ and these cones. If $\beta=0$, then we
define
 $\mathcal A_0(V_1,\dots, V_n)$ as the space of functions
  analytic on the same domain and satisfying
 \begin{equation}
  |\mathbf u(\zeta)|\leq
 C_{R,W_1,\dots,W_n}\prod_{j=1}^n |\I \zeta_j|^{-N}\qquad
   \I \zeta_j\in W_j,\quad |\zeta_j|\leq R,\quad j=1,\dots,n,
 \label{10.2}
 \end{equation}
 where $N$  depends on  $R>0$ and on $W_j$.
 Clearly, these spaces are algebras under pointwise
 multiplication.

 \begin{theorem} The Laplace transformation
  $\mathcal L:\,v\to (v, e^{iz\zeta})$ is an isomorphism
  of the space $S^{\prime\,\beta}(V_1^*,\dots, V_n^*)$,
  $0\le \beta<1$, onto the algebra $\mathcal A_\beta(V_1,\dots, V_n)$.
  The analytic function $(\mathcal
Lv)(\zeta)$ tends to the Fourier transform $\tilde v$ of
 $v$ in the strong topology of $S'_\beta(\oR^{dn})$ as
$\I \zeta\to 0$ inside a fixed cone $W_1\times\dots\times W_n$,
where $W_j\Subset V_j$.
\end{theorem}

\begin{proof} Since $S^{\prime\,0}(V_1^*,\dots, V_n^*)\subset
S^{\prime\,0}(V^*)$, we can use Theorem~4 in~\cite{FS} for
$\beta>0$ and Theorem~2 in~\cite{S} for $\beta=0$. Their
statements are identical to that of Theorem~12 for $n=1$. In
particular, they show that every functional belonging to
$S^{\prime\,\beta}(V^*)$ has a Laplace transform, which is
analytic in $T^V$ and whose boundary value is  $\tilde v$. The
bounds~\eqref{10.1} and \eqref{10.2} are stronger than the bounds
in~\cite{FS,S}, which hold for an arbitrary element of
$S^{\prime\,\beta}(V^*)$. However, the multiplicative
property~\eqref{m} enables us to derive them in the same way,
starting from the estimate
\begin{equation}
|\mathcal Lv(\zeta)| =|(v,e^{iz\zeta})|\le
\|v\|'_{U,B,N}\|e^{iz\zeta}\|'_{U,B,N},
 \label{10.3}
 \end{equation}
 where we use the norms~\eqref{2*} and their  dual norms.
 Here $B$ can be taken arbitrarily large,
 $U=U_1\times\dots\times U_n$, where $U_j$ are any cones with open projections
 such that $V_j^*\Subset U_j$, and  $N$ generally depends on
 $B$ and $U$. We choose the cones $U_j$ and   auxiliary cones $U'_j$
 so that $V^*_j\Subset U_j\Subset U'_j\Subset \Int W^*_j$,
 where $\Int W^*_j$ is the interior of $W^*_j$. This is possible because
 $W_j\Subset V_j$ implies that $V^*_j\Subset \Int W^*_j$.
 Let $\beta=0$ and $\zeta=\xi+i\eta$. Then
   \begin{equation}
   \|e^{iz\zeta}\|'_{U,B,N} = \sup_{x,y}\exp\left\{
   -x\eta-y\xi+N\ln\left(1+|x|\right) - B d(x,U) - B|y|\right\}.
     \label{10.4}
   \end{equation}
This exponential is factorizable, and each factor can be estimated
in the same manner. Namely, assuming that $|\xi_j|\le R<B$, we can
omit  terms that depend on
     $y_j$.  If $x_j\notin U'_j$, then
   $d(x_j,U_j)>\theta |x_j|$ with some $\theta>0$,  and
the expression in the exponent is dominated by a constant for
    $|\eta_j|\le R<\theta B$,. If $x_j\in
U'_j$, then the inclusion $U'_j\Subset\Int W^*_j$ implies that
there is a $\theta'>0$ such that $x_j\eta_j\ge \theta'|x_j|
|\eta_j|$ for all $x_j\in U'_j$ and $\eta_j\in W_j$. Substituting
this inequality in \eqref{10.4}, dropping the negligible term
$d(x_j,U_j)$, and locating the extremum,  we obtain~\eqref{10.2}
with some constant $C_{R,W_1,\dots, W_n}$ proportional to
$\|v\|'_{U,B,N}$. The  case $\beta>0$ is treated in the same way,
with obvious changes in computation.

The nontrivial part of  Theorem 12 states that any function
belonging to the algebra $\mathcal A_\beta(V_1,\dots, V_n)$ is the
Laplace transform of an element in $S^{\prime\,\beta}(V_1^*,\dots,
V_n^*)$.  Let $\mathbf u$ be a function with property~\eqref{10.2}
and let $u$ be its boundary value, which exists in the  Schwartz
space  $\mathcal D'(\oR^{dn})=S'_0(\oR^{dn})$ of distributions by
Theorem~3.1.15 of~\cite{H1}. By Theorem~4 in~\cite{FS}, the
stronger condition~\eqref{10.1} implies that the distribution $u$
belongs to $S'_\beta(\oR^{dn})$. Restricting it to  test functions
of the form $g_1\otimes\dots\otimes g_n$, where $g_j\in
S_\beta(\oR^d)$, and using the same theorem for $\beta>0$ and
Theorem~2 in~\cite{S} for $\beta=0$, we conclude that the
multilinear form determined by the inverse Fourier transform of
$u$ is carried by the cone $V_1^*\times\dots\times V_n^*$. An
application of Theorem~11 completes the proof.
\end{proof}
We also note that every cone  $V$ has the same dual cone as  its
convex hull $\ch V$. Hence Theorem~12 implies that
 $\mathcal A_\beta(V_1,\dots, V_2)=
   \mathcal A_\beta(\ch V_1,\dots,\ch V_n)$.

\bigskip

\begin{center}
{\bf Appendix 1.}
\end{center}
\begin{proof}[Proof of Lemma 1] We first show that for any $\sigma>2$
there is a sequence of functions $\varphi_n\in S^0(\oR)$,
$n=0,1,2,\dots$, such that
 $$
|\varphi_n(z)| \le A_N(1+|x|)^{-N}e^{\sigma|y|},
 \eqno (\mathrm a1)
$$
$$
\ln|\varphi_n(iy)|\ge |y|,
 \eqno (\mathrm a2)
$$
$$
 \ln|\varphi_n(z)| \le \sigma|y|-n\ln^+(|x|/n)+A,
  \eqno (\mathrm a3)
$$
where $\ln^+r=\max (0,\ln r)$ and the constants  $A_N$ and $A$ are
independent of $n$. (We take $n\ln^+(|x|/n)=0$ when $n=0$.)

Such a sequence can be constructed by an iterative procedure used
in the theory of quasi-analytic classes and described, for
example, in~ \cite{H1}, \S~1.3. Let $a_0\ge a_1\ge\dots $ be a
sequence of positive numbers. Let $H_a(t)=a^{-1}$ for $-a/2<t<a/2$
and $H_a(t)=0$ outside this range. We set
$$
\omega_n=H_{a_0}*\dots *H_{a_n}.
 \eqno (\mathrm a4)
$$
Clearly, $\omega_n$ is an even nonnegative function supported in
the interval $|t|\le (a_0+\dots+a_n)/2$. The integral of this
function equals 1 because $\int(u*v)\d t=\int\! u\,\d t\int\!
v\,\d t$. Using the relation
$$
(u*H_a)'(t)=\frac{u(t+a/2)-u(t-a/2)}{a},
$$
where $u$ is  assumed to be a continuous function, we see that
$\omega_n^{(k)}$ can be written as a sum of $2^k$ terms, each of
which  is a shift of the function  $H_{a_k}*\dots
*H_{a_n}/(a_0\dots a_{k-1})$. Taking the inequality $|u*v|\le\sup
|u|\int\!|v|\d t$ into account, we obtain
$$
  |\omega_n^{(k)}| \leq  \frac{2^k}{a_0\dots a_k},\quad 0\leq k\leq n.
  \eqno (\mathrm a5)
  $$
 We note that $\omega_n\in C^{n-1}_0$, and while the higher derivative
  $\omega_n^{(n)}$ is only piecewise continuous, the estimate
  (a5) holds in this case as well.
We set $a_0=2$, $a_1=\dots=a_n=2/n$. Then
$$
  |\omega_n^{(k)}| \leq  \frac{1}{2}\,n^k, \quad k\leq n,\quad
  \int\omega_n\,\d t=1,\quad \supp \omega_n\subset [-2,2].
    $$
Let us consider the convolution $\psi_n=\omega_n*\omega$, where
$\omega\in C^\infty_0$ is a smooth nonnegative even function
supported in $[-\delta, \delta]$ and having integral 1. The
Laplace transform of $\psi_n$ is estimated as follows:
$$
|x^k\tilde \psi_n(z)|\leq \int\limits_{-2-\delta}^{2+\delta}\left|
e^{izt} \psi_n^{(k)}(t)\right|\,\d t\leq \begin{cases}
C_k\,e^{(2+\delta) |y|}& \text{for all $k$;}\\
\frac{1}{2}n^k\,e^{(2+\delta) |y|} &\text{for  $k\le n$.}
\end{cases}
\eqno (\mathrm a6)
 $$
On the other hand, $\int\limits_{|t|>1-\delta}\psi_n(t)\,\d
t\ge\delta$ because $\psi_n\le 1/2$ and $\int\psi_n(t)\,\d t=1$.
Since the function $\psi_n$ is  nonnegative and even, we get
   $$
   |\tilde \psi_n(iy)|=\int\, e^{-yt}\, \psi_n(t)\,\d t\geq
   \int_{t>1-\delta}e^{|y|t} \psi_n(t){\rm d}t\geq
{\delta\over2}e^{(1-\delta)|y|}.
 \eqno (\mathrm a7)
 $$
Hence the sequence $\varphi_n(z)=(2/\delta){\tilde \psi_n}
(z/(1-\delta))$ possesses all the required properties (a1)--(a3),
if $\delta$ is chosen so that $(2+\delta)/(1-\delta)<\sigma$.

This sequence is the main tool for proving Lemma~1. Without loss
of generality, we can assume that $b=1$. Let us introduce the
auxiliary function
$$
 H(\xi)=\sup_y\{\ln|\eta(\xi+iy)|-|y|\},
\eqno (\mathrm a8)
 $$
By~\eqref{11}, it satisfies the inequality
$$
 H(\xi)\le \ln C_N-N\ln(1+|\xi|)+d(\xi,V).
 \eqno (\mathrm a9)
 $$
 We first consider the simplest one-dimensional case, when
$V=\oR_-$ and $d(\xi,V)=\vartheta(\xi)\,|\xi|$, where
$\vartheta(x)$ is the Heaviside step function. Let
$\Phi_n(z)=\ln|\varphi_n(ez)|$. The function $\Phi_n$ is
subharmonic according to~\cite{V}, \S~II.9.12. As a candidate for
the desired function $\varrho$, we take the upper envelope of the
family $\Phi_n(z-\xi)+H(\xi)$, allowing the index $n$  to
 depend on the point $\xi\in\oR$.  The functions in
this family are locally uniformly bounded from above and hence
their upper envelope is also subharmonic (see~\cite{V},
\S~II.9.6). Moreover, it obviously dominates $\ln|\eta(z)|$
because relations (a2) and (a8) imply that
$$
\sup_\xi\{\Phi_n(z-\xi)+H(\xi)\}\ge \Phi_n(iy)+H(x)\ge
\ln|\eta(z)|.
 \eqno (\mathrm a10)
 $$
 We claim that the second inequality in~\eqref{12} is ensured
 by an appropriate choice of $n(\xi)$.
 If $\xi<0$, then $d(\xi,\oR_-)=0$ and we can simply set
$n(\xi)=0$, because the property (a1) and the elementary
inequality
$$
-\ln(1+|x-\xi|)- \ln(1+|\xi|)\le- \ln(1+|x|),
 \eqno (\mathrm a11)
 $$
yield that
$$
\sup_{\xi<0}\{\Phi_0(z-\xi)+H(\xi)\}\le A'_N+e\sigma|y|
-N\ln(1+|x|).
 $$
In view of  (a11), we also have the estimate
$$
\varkappa\,\Phi_0(z-\xi)-N\ln(1+|\xi|)\le A''_N+\varkappa
e\sigma|y| -N\ln(1+|x|), \eqno (\mathrm a12)
 $$
 which holds for any $\varkappa>0$ and  all $\xi$ and shows
  that  difficulties emerge only from the linear growth
 of the term  $d(\xi, V)$ in (a9).

 Suppose that
$\xi\ge 0$ and hence $d(\xi,\oR_-)=\xi$. Suppose also that
$e|x-\xi|> n$. Then
$$
n\ln\frac{e|x-\xi|}{n} +e d(x,\oR_-)\ge n\ln\frac{e\xi}{n}.
 \eqno (\mathrm a13)
 $$
This is obvious for $|x-\xi|>\xi$. When $|x-\xi|\le\xi$, it
suffices to use the inequality $\vartheta(x)\,x\ge\xi-|x-\xi|$ and
note that the function  $n\ln(\lambda/n)-\lambda$ is monotone
decreasing in  $\lambda\in [n,e\xi]$. Combining (a3) and (a13), we
get
$$
\Phi_n(z-\xi)+\xi \le A+\sigma e
|y|+ed(x,\oR_-)-n\ln\frac{e\xi}{n}+\xi.
 $$
 We take  $n(\xi)$ to be the integer part of  $\xi$.
 Then $n\ln(e\xi/n)\ge n>\xi-1$ and
$$
\Phi_{n(\xi)}(z-\xi)+\xi \le  A'+\sigma e |y|+ed(x,\oR_-).
 \eqno (\mathrm a14)
$$
 An analogous inequality holds for   $e|x-\xi|\le n$,
when  $\ln^+(e|x-\xi|/n)$ vanishes. Indeed, in that case $\xi\le
\vartheta(x)\,|x|+|x-\xi|\le \vartheta(x)\,|x|+\xi/e$, and hence
$\xi\le ed(x,\oR_-)$. Thus the inequality (a14) (with an
appropriate constant on the right-hand side) holds for all $\xi\ge
0$. Combining this with estimate (a12) and setting
$\varkappa=B/(e\sigma)-1$ in this estimate, we conclude that the
upper envelope
$$
\varrho(z)=\varlimsup_{z'\to
z}\sup_\xi\{\varkappa\Phi_0(z'-\xi)+\Phi_{n(\xi)}(z'-\xi)+H(\xi)\}
\eqno (\mathrm a15)
 $$
  satisfies all the requirements~\footnote{Taking the upper limit ensures
the upper semicontinuity of the resulting function and enters into
the definition~\cite{V} of upper envelope.}.

In the general case of several variables and an arbitrary open
cone $V\subset\oR^d$, we set
$\Phi_n(z)=\sum_{j=1}^d\ln|\varphi_n(e\sqrt{d}\,z_j)|$. Clearly,
the inequality  (a10) holds. The estimate (a12) is replaced by
$$
 \varkappa\,\Phi_0(z-\xi)-N\ln(1+|\xi|)\le
A'''_N+\varkappa e d\sigma|y| -N\ln(1+|x|),
 \eqno (\mathrm a16)
 $$
 because $\sum_{j=1}^d|y_j|\le\sqrt{d}\,|y|$.
For $x\notin V$ and $e|x-\xi|> n$, we have
$$
\sum_{j=1}^d n\ln\left(\frac{e\sqrt{d}}{n}|x_j-\xi_j|\right) +e
d(x,V)\ge n\ln\left(\frac{e}{n}d(\xi, V)\right).
 \notag
 $$
 To prove this, it suffices to use the formulae
$$
\sum_{j=1}^d \ln^+|x_j|\ge\ln^+\frac{|x|}{\sqrt{d}},\qquad  d(\xi,
V)=\inf_{\xi'\in V}|\xi'-\xi|\le d(x,V)+|x-\xi|.
$$
This time we take $n(\xi)$ to be the integral part of $d(\xi,V)$.
Then  (a14) is replaced by the inequality
$$
\Phi_{n(\xi)}(z-\xi)+d(\xi, V) \le  A''+\sigma e d |y|+ed(x,V),
 $$
 which holds for all $x$. Combining this inequality with  (a16),
 we conclude that the  conditions~\eqref{12} are fulfilled for the
  plurisubharmonic function defined
by (a15) with $\xi\in \oR^d$ and $\varkappa=B/(e\sigma d)-1$.
Lemma~1 is proved.
\end{proof}

\begin{center}
 {\bf Appendix 2.}
\end{center}

\medskip

Suppose that $1<\alpha'<\alpha$. We now return to (a4) and set
$a_0=a_1=1$ and $a_k=(k-1)^{k-1)\alpha'}/k^{k\alpha'}$ for $k>1$.
 The series $\sum a_k$  converges because $(k-1)^{k-1}/k^k\le 1/k$.
 By Theorem~1.3.5  of \cite{H1}, the corresponding sequence (a4)
tends to a smooth nonnegative even function
 as $n\to\infty$. This function is
compactly supported, and its $k$\!~th derivative  is bounded by
$2^kk^{k\alpha'}$. By a scaling transformation it can be converted
into a function $\gamma$ such that
$$
  |\gamma^{(k)}| \leq  C_\epsilon\epsilon^k k^{k\alpha},\qquad
     \supp \gamma\subset [-1/2,1/2],
    $$
where  $\epsilon$ is arbitrarily small. Let $g_\alpha(x)=\tilde
\gamma^2(x)$, where $\tilde \gamma$ is the Laplace transform of
$\gamma$. Then
   $$
    |g_\alpha(x+iy)|\leq  e^{-2|x|^{1/\alpha}+|y|},\quad
  g_\alpha(x)\geq 0,\quad g_\alpha(0)> 0.
  \eqno (\mathrm a17)
   $$

   \begin{proof}[Proof of Theorem 9] We first consider the special case
   of the closed half-plane
  $\oR\times \bar \oR_-$ in $\oR^2$ (the general case can easily
  be  reduced to this one). Using $g_\alpha$, we can construct functions
 $f_1\in S^0(\oR)$ and $f_2\in  S^0(\oR_-)$ such that
   $f_1\otimes f_2\in S^0(\oR, \bar \oR_-)=S^0(\oR\times
   \oR_-)$ and $f_1\otimes f_2\not\in S^\beta(\oR\times\bar \oR_-)$
    for any $\beta\in [0,1)$.   We set
     $$ f_1(x)=\int
  e^{-|\xi|^{1/\alpha}}g_\alpha(x-\xi)\,{\rm d}\xi.
   $$
   This convolution can be analytically continued to whole of $\oC$
   and belongs to    $S^0_\alpha(\oR)$.  Indeed,
using the triangle inequality for the metric $|x-\xi|^{1/\alpha}$,
we obtain
  $$
  |f_1(x+iy)|\leq \int
  e^{-|x-\xi|^{1/\alpha}-2|\xi|^{1/\alpha}+|y|}\,{\rm d}\xi\leq
  C\, e^{-|x|^{1/\alpha}+|y|},
  $$
     where $C=\int e^{-|\xi|^{1/\alpha}}\,{\rm d}\xi$.
  In addition, we have the  lower estimate
  $$
  f_1(x)\geq \int_{-1}^{+1}e^{-|x-\xi|^{1/\alpha}}g_\alpha(\xi)\,{\rm d}\xi
  \geq e^{-(|x|+1)^{1/\alpha}}\int_{-1}^{+1}g_\alpha(\xi)\,{\rm d}\xi\geq
   c\,e^{-|x|^{1/\alpha}}.
  \eqno (\mathrm a18)
  $$
  Furthermore, let $1<\alpha'<\alpha$ and let
  $$
     f_2(x)=\int_0^\infty e^{\xi^{1/\alpha'}}g_{\alpha'}(x-\xi)\,{\rm d}\xi.
  $$
  In an analogous way, it is easy to verify that
   $$
  |f_2(x+iy)| \leq  C'\, e^{|x|^{1/\alpha'}+|y|}.
  \eqno (\mathrm a19)
  $$
 For $x> 1$, we have the estimate
  $$
  f_2(x)= \int_{-\infty}^{x}e^{(x-\xi)^{1/\alpha'}}g_{\alpha'}(\xi)\,{\rm
  d}\xi \geq e^{(x-1)^{1/\alpha'}}\int_{-1}^{+1}g_{\alpha'}(\xi)\,{\rm
  d}\xi\geq c'\,e^{x^{1/\alpha'}}.
  \eqno (\mathrm a20)
  $$
  If $x< 0$ and $\xi>0$, then $|x-\xi|=|x|+|\xi|$. Using the
  inequality $2(|x|+|\xi|)^{1/\alpha'}\geq
  |2x|^{1/\alpha'}+|2\xi|^{1/\alpha'}$, we obtain
  $$
  |f_2(x+iy)|\leq
  C^{\prime\prime}\, e^{-|2x|^{1/\alpha'}+|y|},  \qquad  x\in  \oR_-.
  \eqno (\mathrm a21)
  $$
  The estimates (à19) and (à21) imply that $f_2\in S^0(\oR_-)$
  and,  therefore,
  $$
  (f_1\otimes f_2)(x_1, x_2)=f_1(x_1)f_2(x_2)\in S^0(\oR\times
  \oR_-).
  $$
  On the other hand, the lower estimates (a18) and (a20) show that the function
   $f_1\otimes f_2$ increases to infinity along any real ray in
   the half-plane $x_2>0$. Hence it does not belong to any of the
   $S^\beta(U)$, where $U\Supset \oR\times\bar \oR_-$. Setting
   $$
   (v,f)=\int_1^\infty f(x_1, x_1^{\alpha'/\alpha})\,{\rm d}x_1,
  $$
   we obtain a simple example of a functional in
     $S^{\prime\, 0}(\oR\times \bar \oR_-)$ which does not belong to
      $S^{\prime\, 0}(\oR, \bar \oR_-)$. The function
   $f_1\otimes f_2$ constructed above is bounded below by a positive constant
    on the path of integration  $x_2=x_1^{\alpha'/\alpha}$,
    $x_1>1$. Therefore  $v$  has no continuous extension to
  $S^0(\oR \times \oR_-)$, nor  to any $S^\beta(\oR \times \oR_-)$,
  $\beta>0$. Indeed,  $f_1\otimes f_2$  can be approximated
  (in the topology of $S^0(\oR \times \oR_-)$) by positive functions
  $f_\nu=f_1\otimes (g_\nu\,f_2)\in S^0_\alpha(\oR^2)$, where
  $g_\nu(x_2)=g_{\alpha^{\prime\prime}}(x_2/\nu)$,
  $\alpha^{\prime\prime}<\alpha'$, and the normalization condition
  $g_{\alpha^{\prime\prime}}(0)=1$ is assumed. Clearly we have
  $(v,f_\nu)\to \infty$ as  $\nu\to \infty$, so there is no continuous extension.

  Now, let $K_1$ and $K_2$ be  closed cones in $\oR^{d_1}$ and
   $\oR^{d_2}$, where $d_1\ge 1$ and $d_2\ge 1$.
Suppose  that $K_1\neq \{0\}$ and $K_2\neq \oR^{d_2}$. We claim
that $S^\beta(K_1\times  K_2)$ does not contain $S^0(K_1, K_2)$
and does not even contain the smaller  space $S^0(\oR^{d_1},
K_2)$. Indeed, assume that the first basis vector  $e^1_1$
  in  $\oR^{d_1}$ belongs to $K_1$ and the basis vector $e^1_2$ in
  $\oR^{d_2}$ does not belong to $K_2$.  Let $h_1$ be a function
  in $S^0_\alpha(\oR^{d_1-1})$ such that $h_1(0)\ne 0$ and replace
  $f_1$  by  $f_1\otimes  h_1$ in  the above construction. Clearly,  $f_1\otimes h_1\in S^0_\alpha(\oR^{d_1})$.
  We also replace $f_2$ by
  $f_2\otimes h_2$, where $h_2\in S^0_{\alpha^{\prime\prime}}(\oR^{d_2-1})$
  and $h_2(0)\ne 0$. It is easy to see that
     $f_2\otimes h_2\in S^0_{\alpha'}(U)$, where  $U$
     is an  open cone  in  $\oR^{d_2}$  defined by the inequality
  $(1+\theta)x^1_2<|x_2|$.
      Clearly, $K_2\setminus\{0\}$ is contained in this cone
       if $\theta>0$ is small enough.
       Therefore, the function $f_1\otimes
  h_1\otimes f_2\otimes h_2$ belongs to $S^0_\alpha(\oR^{d_1}\times U)$,
  while none of the spaces  $S^\beta(K_1\times  K_2)$,
$0\le\beta<1$, contains this function, as is evident from (a18)
and (a20).  In complete analogy with what was done above, we
define a functional $v$ by  integrating test functions
  along the curve $x^1_2=(x^1_1)^{\alpha'/\alpha}$, $x_1>1$, in the
  plane $\{e^1_1,  e^1_2\}$.  This functional is carried by the ray
$\{\lambda e^1_1\mid \lambda\geq 0\}$ lying on the boundary of
$K_1\times K_2$, but $v$ does not belong to
$S^{\prime\,\beta}(K_1, K_2)$, nor even to
$S^{\prime\,\beta}(\oR^{d_1}, K_2)$. This completes the proof.
\end{proof}

   \begin{acknowledgements} This work was supported in part by
the Russian Foundation for Basic Research (grant no.~05-01-01049)
and by the Programme for the Support of Leading Scientific Schools
(grant no.~LSS-4401.2006.2).
\end{acknowledgements}

\end{document}